\documentclass[english,a4paper,12pt]{article}
\usepackage{rotating} 
\usepackage{color}
\usepackage{doc}
\usepackage{graphicx}
\usepackage{caption}
\usepackage{subcaption}
\usepackage{amsmath}
\usepackage{amsfonts}
\usepackage{amssymb}
\usepackage{babel}
\usepackage{latexsym}
\usepackage[utf8]{inputenc}

\newcommand{\sgn}{\mathop{\mathrm{sgn}}}
\title{The foundations of fractional Mellin transform analysis\thanks{This article  is a far-reaching  extension of 
the paper 
" Mellin analysis and exponential sampling. Part I: Mellin fractional integrals
" ( jointly with C. Bardaro and P.L. Butzer, see \cite{BBM} ),
presented by Ilaria Mantellini at  "SAmpta2013", held at Jacobs University, Bremen, on 1-5 July, 2013, and was conducted by 
Goetz Pfander, Peter Oswald, Peter Massopust and Holger Rauhut.
The meeting marked two decades of  SAmpta Workshops
[held at  Riga (1995), Aveiro (1997),  Loen (1999), Tampa (2001),
Strobl (2003), Samsun (2005), Marseille 
(2007), Singapore (2009) ] and the 
Commemoration of the  85-th Birthday of P. L. Butzer.
}}
\author{Carlo Bardaro \thanks{
Department of Mathematics and Computer Sciences, University of Perugia,
via Vanvitelli 1, I-06123 Perugia, Italy, e-mail:
carlo.bardaro@unipg.it; corresponding author} \and Paul L. Butzer \thanks{Lehrstuhl A fur Mathematik RWTH-Aachen
Templergraben 55, D-52062 Aachen, Germany,
e-mail:
butzer@rwth-aachen.de}
\and  Ilaria Mantellini \thanks{Department of Mathematics and Computer Sciences, University of Perugia,
via Vanvitelli 1, I-06123 Perugia, Italy, e-mail:
ilaria.mantellini@unipg.it }}

\renewcommand{\indent}{\hskip18pt}
\date{{\it{In Memory of Rashid Gamid-oglu Mamedov, a pioneer in Mellin Analysis}}}

\begin{document}
\maketitle
\noindent
{\bf Abstract:} {\small In this article we study the basic theoretical properties of Mellin-type fractional integrals, known as generalizations of the Hadamard-type fractional integrals. 
We give a new approach and version, specifying their semigroup property, their domain and range. Moreover we introduce a notion of strong fractional Mellin derivatives and we study the connections with the pointwise fractional Mellin derivative, which is defined by means of Hadamard-type fractional integrals. One of the main results is a fractional version of the fundamental theorem of differential and integral calculus in the Mellin frame. In fact, in this article it will be shown that the very foundations of Mellin transform theory and the corresponding analysis are quite different to those of the Fourier transform, alone since even in the simplest non-fractional case the integral operator (i.e. the anti-differentiation operator) applied to a function $f$ will turn out to be the $\int_0^x f(u)du/u$ with derivative $(xd/dx)f(x).$  Thus the fundamental theorem in the Mellin sense is valid in this form, one which stands apart from the classical Newtonian integral and derivative.
Among the applications two fractional order partial differential equations are studied.}
\vskip0,2cm
\noindent
\vskip0,1cm
\noindent

\section{\bf Introduction}

The theory of Mellin transforms as well as Mellin approximation theory  was introduced by R.G. Mamedov in his treatise \cite{MA}, which includes also previous results in this subject obtained in collaboration with G.N. Orudzhev (see \cite{MO0, MO1, MO2}). In his review Professor H.J. Glaeske (MR1235339--94:44003) writes: {\em This book deals with the theory of The Mellin transform and its applications to approximation theory based on results of the school of I.M. Dzhrbashyan and the methods of the school of P.L. Butzer on Fourier Analysis and approximation}. Somewhat later Mellin transform theory was presented in a systematic form, fully independently of Fourier analysis, by Butzer and Jansche in their papers \cite{BJ}, \cite{BJ1}. Further important developments were then given in \cite{BJ2}, and later on in the present line of research in \cite{BM,BM1,BM2,BM3,BM4, BM5, BM6, MAN, AV1,AV2}.

In the papers \cite{BKT,BKT1,BKT2,BKT3,BKT4} a broad study of fractional Mellin analysis was developed in which the so-called Hadamard- type integrals, which represent the appropriate extensions of the classical Riemann-Liouville and Weyl fractional integrals, are considered (see also the book \cite{KST}). These integrals are also connected with the theory of moment operators (see \cite{BM6},\cite{BM5}, \cite{BCS}).
The purpose of this article is not only a continuation of these topics but also to present a new, almost independent approach, one starting from the very foundations.
As remarked in \cite{BKT}, in terms of Mellin analysis, the natural operator of fractional integration is not the classical Riemann-Liouville fractional integral of order $\alpha >0$ on $\mathbb{R}^+$, namely (see \cite{SKM}, \cite{MSM}, \cite{BW}, \cite{BW1})
\begin{equation}
(I^\alpha_{0+}f)(x) = \frac{1}{\Gamma (\alpha)}\int_0^x (x-u)^{\alpha -1}f(u)du~~(x>0)
\end{equation}
but the Hadamard fractional integral, introduced essentially by Hadamard in \cite{HAD},
\begin{equation}
(J^\alpha_{0+} f)(x) = \frac{1}{\Gamma(\alpha)} \int_0^x  \bigg( \log \frac{x}{u} \bigg)^{\alpha-1} f(u) \frac{du}{u}~~(x >0).
\end{equation}
Thus not $\int_0^xf(u)du$ is the natural operator of integration (anti-differentiation) in the Mellin setting, but $\int_0^x f(u)\frac{du}{u}$ (the case $\alpha =1$).
It is often said that a study of Mellin transforms as an independent discipline is fully superfluous since one supposedly can reduce its theorems and results to the corresponding ones of Fourier analysis by a simple change of variables and functions. It may be possible to reduce a formula by such a change of operations but not the precise hypotheses under which a formula is valid. But alone since (1) is not the natural operator of integration in the Mellin frame but that the Hadamard fractional integral (2), (which is a compact form of the iterated integral (6) (see section 4) ) will turn out to be the operator of integration, thus anti-differentiation to the operator of differentiation $D_{0+,0}f$ in (4) (see below)--in the sense that the fundamental theorem of the differential and integral calculus must be valid in the Mellin frame--makes the change of operation argument fully obsolete. This will become evident as we proceed along, especially in Theorems 3-4, and Theorems 6-12 below. Thus the very foundations to Mellin analysis are quite different to those of classical Fourier analysis.

For the development of the theory, it will be important to consider the following generalization of the  fractional integral, known as the Hadamard-type fractional integrals, for $\mu \in \mathbb{R},$ namely (see \cite{BKT,BKT1,BKT2,BKT3,BKT4, KST})
\begin{equation}
(J^\alpha_{0+, \mu} f)(x) = \frac{1}{\Gamma(\alpha)} \int_0^x \bigg( \frac{u}{x}\bigg)^\mu \bigg( \log \frac{x}{u} \bigg)^{\alpha-1} f(u) \frac{du}{u}~~(x>0)
\end{equation}
for functions belonging to the space $X_c$ of all measurable complex-valued functions defined on $\mathbb{R}^+,$ such that $(\cdot)^{c-1}f(\cdot) \in L^1(\mathbb{R}^+).$
As regards the classical Hadamard fractional integrals and derivatives, some introductory material about fractional calculus in the Mellin setting was already treated in \cite{MA} and \cite{SKM}.

In Section 2 we recall some basic tools and notations of Mellin analysis, namely the Mellin transform, along with its fundamental properties, the notion of the basic Mellin translation operator, which is now defined via a dilation operator instead of the usual traslation (see \cite{BJ}. For other classical references see \cite{BRY}, \cite{GPS}, \cite{PBM}, \cite{SZ}, \cite{ZE}, \cite{ZA}).

In Section 3 we will introduce and study a notion of a strong fractional derivative in the spaces $X_c,$ which represents an extension of the classical strong derivative of Fourier analysis in $L^p-$spaces  (see \cite{BN}). The present notion is inspired by an analogous construction given in \cite{WE}, \cite{BW} for the Riemann-Liouville fractional derivatives in a strong sense. This method is based on the introduction of certain fractional differences, which make use of the classical translation operator. Another important fact is that fractional differences are now defined by an infinite series. Our definition here, follows this approach, using the Mellin translation operator. Our definition reproduces the Mellin differences of integral order, as given in \cite{BJ2}, in which we have a finite sum.

It should be noted that a different approach for spaces $X_0$ was introduced in \cite{MA}, pages 175-176,  starting with the incremental ratios of the integral (2).
A relevant part of the present paper (Section 4) deals with the pointwise fractional derivative of order $\alpha >0,$ known as the "Hadamard-type fractional derivative"
in the local spaces $X_{c,loc},$ and with its links with the strong derivatives.
This notion originates from the analogous concept of Riemann-Liouville theory, and was introduced in \cite{BKT1} using the Hadamard-type fractional integrals. It read as follows
\begin{equation}
(D^\alpha_{0+,\mu} f)(x) = x^{-\mu} \delta^m x^\mu (J^{m- \alpha}_{0+, \mu} f)(x),
\end{equation}
where $m = [\alpha] +1$ and $\delta:=(x \frac{d}{dx})$ is the Mellin differential operator $(\delta f)(x) = x f'(x),$ provided $f'(x)$ exists. For $\mu = 0$ we have the so called Hadamard fractional derivative, treated also in \cite{MA}, \cite{SKM}.
Note that the above definition reproduces exactly the Mellin derivatives $\Theta_c^kf$ of integral order when $\alpha = k \in \mathbb{N}.$ Thus
$D^\alpha_{0+,\mu} f$ represents the natural fractional version of the differential operator $\Theta_c^k,$ in the same way that the Riemann-Liouville fractional
derivative is the natural extension of the usual derivative.
Paper \cite{KI}, gives some sufficient conditions for the existence of the pointwise fractional derivative for functions defined in bounded intervals $I \subset \mathbb{R}^+,$
involving spaces of absolutely continuous functions in $I.$

Since the definition of the pointwise fractional derivatives is based on a Hadamard-type integral, it is important to study in depth the domain and the range of these integral operators.
As far as we are aware this was not sufficiently developed in the literature so far. Here we define the domain of the operator (3) as the subspace of all functions such that
the integral exists as a Lebesgue integral. A basic result in this respect is the semigroup property of $J^\alpha_{0+,c}.$ This was first studied in \cite{MA} and \cite{SKM} for the Hadamard integrals (2) and then developed for the integrals (3) in \cite{BKT2} and  \cite{KI} (see also the recent books \cite{KST}, \cite{BDST}). However, the above property was studied only for functions belonging to suitable subspaces of the domain, namely the space $X^p_{c}$ of all the functions
$f: \mathbb{R}^+ \rightarrow \mathbb{C}$ such that $(\cdot)^{c-1}f(\cdot) \in L^p(\mathbb{R}^+),$ or for $L^p(a,b)$ where $0< a < b < \infty.$

Here we prove the semigroup property in a more general form, using minimal assumptions. This extension enables us to deduce the following chain of inclusions for the domains of the operators $J^\alpha_{0+,c}.$
$$Dom J^\beta_{0+,c} \subset X_{c,loc} = Dom J^1_{0+, c} \subset Dom J^\alpha_{0+,c},$$
for $\alpha < 1 < \beta,$ and all inclusions are strict.

Concerning the range, we show that $J^\alpha_{0+,c}f \in X_{c,loc}$ whenever $f\in Dom J^{\alpha +1}_{0+,c}f$ and in general
$f \in Dom J^\alpha_{0+,c}$ does not imply that $J^\alpha_{0+,c}f \in X_{c,loc}.$

For spaces $X_c$ we have the surprising result that $J^\alpha_{0+,c}f \not \in X_c$ for any nontrivial non-negative function $f \in Dom J^\alpha_{0+,c}.$ This fact gives problems for the evaluation of the Mellin transform of the function $J^\alpha_{0+,c}f.$ In order to avoid this problem, we prove that if
$f \in Dom J^\alpha_{0+,c} \cap \bigcap_{\mu \in [\nu,c]}X_\mu,$ then $J^\alpha_{0+,c}f \in X_\nu$ and so its Mellin transform can be evaluated on the line
$\nu + it,$ with $\nu < c.$

We then apply the theory to deduce one of the main results of this paper, namely the fundamental theorem of the fractional differential and integral calculus in the Mellin frame, here established under sharp assumptions. We consider also some more general formulae, involving different orders of fractional integration and differentiation. Similar results were also given in \cite{KI}, \cite{KST} however in resrticted subspaces (see the remarks in Section 4). In particular, one of the two fundamental formulae is given there under the strong assumption that the functions $f$ belongs to the range of $J^\alpha_{0+, \mu}(X^p_{0+,c}),$ with $\mu > c.$

In Section 5 we prove an equivalence theorem with four equivalent statements, which connects fractional Hadamard-type integrals, strong and pointwise fractional Mellin derivatives and the Mellin transform (see Theorem 8 below).
As far as we know, a fundamental theorem with four equivalent assertions in the form presented
here for the Mellin transform in the fractional case has never been stated explicitly, for the Fourier transform.
As a fundamental theorem in the present sense it was first established for $2\pi-$periodic functions via the finite Fourier transform in \cite{BW}, and for the Chebyshev transform (see e.g. \cite{GPS}, pp. 116-122), in \cite{BS}, \cite{BS2}.
Fractional Chebyshev derivatives were there defined in terms of fractional order differences of the Chebyshev
translation operator, the Chebyshev integral by an associate convolution product.
The next fundamental theorem, after that for Legendre transforms (see e.g. \cite{GPS}, pp.122-131; \cite{BUT}, \cite{SW}), was the one concerned with the Jacobi transform, see e.g. \cite{BSW}. In their inimitable book \cite{GPS}, H.J. Glaeske, A.P. Prudnikov and K.A. Skornik study the Mellin transform and its essential properties (pp. 55-67), not as an independent discipline but by making use of the corresponding properties of the Fourier transform, the reduction being carried out with unusual precision. In other respects their presentation is standard. Thus their integral is the classical one, i.e. $F(x) = \int_0^x f(u)du,$ with Mellin transform $M[F](s) = -s^{-1}M[f](s+1).$  They were not aware of \cite{BJ}. However, their sections on the Chebyshev, Legendre, Gegenbauer and Jacobi transforms make interesting reading and are unorthodox. Here their chief properties are based on the definitions of an associated translation operator for each transform, an approach carried out systematically for the Chebyshev and Legendre transforms in \cite{BS2} and \cite{SW}, which are cited by the three authors. However, they do not continue the process and define the associated derivative concepts in terms of the respective translation operators (probably due to lack of space). This would have led them to the fundamental theorems of the differential and integral calculus in the setting of the respective transforms. Neverthless the material of these sections has never been treated in a book-form as yet. The chapter on Mellin transforms in the unique handbook \cite{ZA}, also written in the classical style, bears the individuaal stamp of the author, A. Zayed.

In Section 6 we describe some special cases of interest in applications, while in Section 7 we apply our theory to two fractional partial differential equations.
The use of Mellin transforms for solving partial differential equations originates from certain boundary value problems in wedge-shaped regions, see e.g. \cite{ZE}, \cite{KN} and, in the fractional frame, was considered by various authors for the study of fractional versions of the diffusion equation (see e.g. \cite{WY}, \cite{SWY}, \cite{HIL}, \cite{KST}). However, the use of Mellin transforms for solving fractional differential equations with Hadamard derivatives is not usual. Also, there are a few contributions dealing with pure Hadamard derivatives (see e.g. \cite{KST}, \cite{BDST}, \cite{KLY}, \cite{QFT}, \cite{GM}).
Most fractional equations, are studied using different types of fractional derivatives, (Riemann-Liouville, Caputo, etc).
Here we apply our theory to an integro-differential equation which can be reduced to a fractional evolution equation, with Hadamard fractional derivative. A similar equation was also considered in \cite{KST} but with the Caputo fractional derivative. Here we give the exact solution of the evolution equation, using just Mellin transforms and the fractional theory developed in this paper. As a second example, we consider a boundary value problem for a fractional diffusion equation, using the same approach. In both the examples the (unique) solution is given in terms of a Mellin convolution operator.

In the very recent book \cite{BDST} numerical methods for solving fractional differential equation are treated, using mainly Caputo and Riemann-Liouville fractional theories.

\section{Preliminaries}
Let $L^1 = L^1(\mathbb{R}^+)$ be the space of all Lebesgue measurable and integrable complex-valued functions defined on $\mathbb{R}^+,$ endowed with the usual norm.

Let us consider the space, for some $c \in \mathbb{R},$
$$ X_c = \{ f: \mathbb{R}^+ \rightarrow \mathbb{C} : f(x) x^{c-1} \in L^1 (\mathbb{R}^+ )\}$$
endowed with the norm
$$\|f\|_{X_c} = \|f(\cdot) (\cdot)^{c-1} \|_{L^1} = \int_0^\infty |f(u)| u^{c-1} du.$$
More generally by $X^p_c$ we denote the space of all functions $f: \mathbb{R}^+ \rightarrow \mathbb{C}$ such that
$(\cdot)^c f(\cdot) \in L^p(\mathbb{R}^+),$ with $1< p < \infty.$
In particular when $c=1/p,$ the space $X^p_c$ coincides with the classical $L^p(\mathbb{R}^+)$ space.

\noindent
For $a,b \in \mathbb{R}$ we define the spaces $X_{(a,b)},~X_{[a,b]}$ by
$$X_{(a,b)} = \bigcap_{c \in ]a,b[}X_c,~~X_{[a,b]} = \bigcap_{c \in[a,b]}X_c$$
and, for every $c$ in $(a,b)$ or $[a,b]$, $\|f\|_{X_c}$ is a norm on them.

Note that, for any $a,b \in \mathbb{R},$ with $a<b,$ if $f \in X_a \cap X_b$, then $f \in X_{[a,b]}$ and moreover
$$\|f\|_{X_c} \leq \|f\|_{X_a} + \|f\|_{X_b},$$
for every $c \in [a,b].$
For these and other results see \cite{BJ}.
\vskip0,2cm
In what follows, we denote by $\chi_A(x)$ the characteristic function of the set $A \subset \mathbb{R}^+.$
\vskip0,2cm

We define for every $f\in X_c$ the Mellin transform $[f]^\wedge_M$ of $f$ by
$$ M[f](s) \equiv [f]^{\wedge}_M (s) = \int_0^\infty u^{s-1} f(u) du $$
where $s=c+ it, t\in \mathbb{R}.$

The notation $M[f(\cdot)](s)$ of the Mellin transform signifies the fact that one of its essential roles is to solve
analytical problems by transforming them into another function space, solve the problem (which should be simpler) in the trasformed state, and then apply a (suitable) Mellin inversion formula to obtain the solution in the original function space.

Basic in this respect are the linearity and boundedness properties, thus
$$M[a f(\cdot) + b g(\cdot)](s) = a M[f(\cdot)] (s)+ b M[g(\cdot)](s)~~(f,g \in X_c,~a,b \in \mathbb{R})$$
$$|M[f(\cdot)] (s)| \leq \|f\|_{X_c}~~~(s = c+it).$$
As a  consequence of the boundedness property, if $(f_n)_n$ is a sequence of functions in $X_c$ convergent in $X_c$ to a function $f,$ then
$M[f_n]$ converges uniformly to $M[f]$ on the line $s=c+it,~t \in \mathbb{R}.$
\vskip0,4cm

We need several operational properties.

The Mellin translation operator $\tau_h^c$, for $h \in \mathbb{R}^+,~c \in \mathbb{R},$~$f:\mathbb{R}^+ \rightarrow \mathbb{C},$ is defined by
$$(\tau_h^c f)(x) := h^c f(hx)~~(x\in \mathbb{R}^+).$$

Setting $\tau_h:= \tau^0_h,$ then
$$(\tau_h^cf)(x) = h^c (\tau_hf)(x),~\|\tau_h^c f\|_{X_c} = \|f\|_{X_c},~(\tau_h^c)^j f(x) = h^{jc} f(h^j x) = (\tau_{h^j}^c f)(x).$$
 Proposition 2 and Lemma 3 in \cite{BJ}, state the following:
\newtheorem{Lemma}{Lemma}
\begin{Lemma}\label{Lemma 1}
The Mellin translation operator $\tau_h^{\overline{c}}: X_c \rightarrow X_c$~for $c,\overline{c} \in \mathbb{R},~h \in \mathbb{R}^+$
is an isomorphism with $(\tau_h^{\overline{c}})^{-1} = \tau_{1/h}^{\overline{c}}$ and
$$\|\tau_h^{\overline{c}}f\|_{X_c} = h^{\overline{c} - c}\|f\|_{X_c}~~(f \in X_c)$$
having the properties
\begin{itemize}
\item[i)] $M[\tau_h^{\overline{c}}f] (s) =h^{\overline{c} - s}M[f](s),$
in particular $M[\tau_h f](s) = h^{-s}M[f](s);$
\item[ii)] $\lim_{h \rightarrow 1} \|\tau_h^{\overline{c}}f - f\|_{X_c} =0.$
\end{itemize}
\end{Lemma}
\vskip0,3cm
When $\overline{c}=0$ Property ii), in case of continuous functions $f$, expresses uniform continuity in the Mellin frame, taking the usual $L^\infty-$norm, i.e.
$$\lim_{h \rightarrow 1}\|\tau_hf - f\|_{\infty} =0.$$
It is equivalent to the so-called log-uniform continuity due to Mamedov (see \cite{MA}, page 7), which may be expressed as follows: a function $f:\mathbb{R}^+ \rightarrow \mathbb{C}$ is
log-uniformly continuous on $\mathbb{R}^+$ if for every $\varepsilon >0$ there exists $\delta_\varepsilon >0$ such that $|f(u) - f(v)| < \varepsilon,$
whenever $|\log u - \log v| < \delta_\varepsilon.$ Indeed the continuity of the operator $\tau_h$ implies that $|f(hx) - f(x)| < \varepsilon,$ for $|h|< \delta_\varepsilon,$ uniformly with respect to $x \in \mathbb{R}^+.$ It should be noted that this notion is different from the usual uniform continuity.
For example, the function $f(u) = \sin u$ is obviously uniformly continuous, but not log-uniformly continuous on $\mathbb{R}^+$, while the function
$g(u) = \sin (\log u)$ is log-uniformly continuous but not uniformly continuous on $\mathbb{R}^+.$ However, the two notions are equivalent on every bounded
interval $[a,b]$ with $a>0.$
\vskip0,5cm
The Mellin convolution product, denoted by $f\ast g$, of two functions $f,g :\mathbb{R}^+ \rightarrow \mathbb{C},$ is defined by
$$(f\ast g)(x) := \int_0^{+\infty} g(\frac{x}{u})f(u)\frac{du}{u} = \int_0^{+\infty}(\tau^c_{1/u}f)(x)g(u)u^c\frac{du}{u}~~~(x \in \mathbb{R}^+)$$
in case the integral exists. It has the properties
\begin{Lemma}\label{convolution}
\begin{itemize}
\item[(i)] If $f,g \in X_c,$ for $c \in \mathbb{R},$ then $f\ast g$ exists (a.e.) on $\mathbb{R}^+,$ it belongs to $X_c$, and
$$\|f \ast g\|_c \leq \|f\|_{X_c} \|g\|_{X_c}.$$
If in addition $x^c f(x)$ is uniformly continuous on $\mathbb{R}^+,$ then $f\ast g$ is continuous on $\mathbb{R}^+.$
\item[(ii)] (Convolution Theorem). If $f,g \in X_c$ and $s=c+it,~t \in \mathbb{R},$ then
$$M[f \ast g](s) = M[f](s) M[g](s).$$
\item[(iii)] (Commutativity and Associativity). The convolution product is commutative and associative, thus for $f_1,f_2, f_3 \in X_c$ there holds true (a.e.)
$$f_1 \ast f_2 = f_2 \ast f_1,~~(f_1 \ast f_2) \ast f_3 = f_1 \ast (f_2 \ast f_3).$$
In particular $X_c$ is a Banach algebra.
\end{itemize}
\end{Lemma}
\section{The strong Mellin fractional differential operator}

Let us denote by $I$ the identity operator over the space of all measurable functions on
$\mathbb{R}^+.$

The Mellin fractional difference of $f \in X_c$ of order $\alpha >0,$ defined by
\begin{eqnarray*}
\Delta_h^{\alpha, c} f(x):= (\tau_h^c - I)^\alpha f(x) = \sum_{j=0}^\infty \left(\begin{array}{c} \alpha\\ j \end{array} \right) (-1)^{\alpha-j} \tau_{h^j}^c f(x)
\end{eqnarray*}
for $h>0$
with
$$\left(\begin{array}{c} \alpha\\ j \end{array} \right) = \frac{\alpha (\alpha -1)\cdots(\alpha -j +1)}{j!},$$
has the following properties
\newtheorem{Proposition}{Proposition}
\begin{Proposition}\label{Proposition 1}
For $f \in X_c$ the difference $\Delta_h^{\alpha,c} f(x)$ exists a.e. for $h >0,$ with
\begin{itemize}
\item[i)] $\|\Delta_h^{\alpha,c} f\|_{X_c} \leq \|f\|_{X_c}\sum_{j=0}^\infty \bigg|\left(\begin{array}{c} \alpha\\ j \end{array} \right)\bigg|$
\item[ii)] $M[\Delta_h^{\alpha,c} f](c+it) = (h^{-it} -1)^\alpha M[f](c+it).$
\item[iii)] The following semigroup property holds for $\alpha, \beta >0,$
$$(\Delta_h^{\alpha,c}\Delta_h^{\beta,c}f)(x) = (\Delta_h^{\alpha + \beta,c}f)(x).$$
\end{itemize}
\end{Proposition}
{\bf Proof}.  At first, we have for $x >0,~h>0$
\begin{eqnarray*}
|\Delta_h^{\alpha,c} f (x)| \leq \frac{1}{x^c}\sum_{j=0}^\infty \bigg | \left(\begin{array}{c} \alpha\\ j \end{array} \right)\bigg| h^{cj}x^c|f(h^jx)|;
\end{eqnarray*}
thus we have to prove the convergence of the latter series.
For this purpose, by integration, we have
\begin{eqnarray*}
&&\int_0^\infty \sum_{j=0}^\infty \bigg | \left(\begin{array}{c} \alpha\\ j \end{array} \right)\bigg|(h^j x)^c |f(h^jx)|\frac{dx}{x} = \sum_{j=0}^\infty \bigg | \left(\begin{array}{c} \alpha\\ j \end{array} \right)\bigg| \int_0^\infty (h^j x)^c |f(h^jx)| \frac{dx}{x} := J.
\end{eqnarray*}
Now, putting in the second integral $h^jx = t$, we have
$$J =
 \|f\|_{X_c}\sum_{j=0}^\infty \bigg | \left(\begin{array}{c} \alpha\\ j \end{array} \right)\bigg|.
$$
Thus, since
$\left(\begin{array}{c} \alpha\\ j \end{array} \right) = {\cal O}(j^{-\alpha - 1}),~~j \rightarrow + \infty,$
we observe that the integral is finite for any $h>0,$ if $f \in X_c,$ and so the integrand is finite almost everywhere. Thus the original series defining the difference, converges almost everywhere.

As to (i), we have
\begin{eqnarray*}
&&\|\Delta_h^{\alpha,c} f\|_{X_c} = \int_0^{\infty} x^{c-1}\bigg|\sum_{j=0}^\infty \left(\begin{array}{c} \alpha\\ j \end{array} \right) (-1)^{\alpha -j} h^{cj} f(h^j x)  \bigg|dx \\
&&\leq \int_0^{\infty} \frac{t^{c-1}}{h^{j(c-1)}} \sum_{j=0}^\infty \bigg|\left(\begin{array}{c} \alpha\\ j \end{array} \right)\bigg|  h^{cj} |f(t)| \frac{dt}{h^j} = \|f\|_{X_c} \sum_{j=0}^\infty \bigg|\left(\begin{array}{c} \alpha\\ j \end{array} \right)\bigg|,
\end{eqnarray*}
and so the assertion.

An alternative proof makes use of Lemma \ref{Lemma 1} in the following way.
The left hand side of i) can be estimated by
$$\|\Delta_h^{\alpha,c} f\|_{X_c} \leq \sum_{j=0}^\infty \bigg|\left(\begin{array}{c} \alpha\\ j \end{array} \right)\bigg| \|\tau_{h^j}^cf\|_{X_c} =  \sum_{j=0}^\infty \bigg|\left(\begin{array}{c} \alpha\\ j \end{array} \right)\bigg|\|f\|_{X_c}$$
which is independent of $h>0.$
As to (ii), the Mellin transform on the left equals, by the linearity property, which uses an integration by series,
$$\sum_{j=0}^\infty \left(\begin{array}{c} \alpha\\ j \end{array} \right) (-1)^{\alpha -j}h^{-itj}[f]^\wedge_M(c+it),$$
which yelds (ii).
Note that the complex number $h^{-it}$ has modulus $1,$ so it lies on in the boundary of the circle of convergence
of the power series which defines the binomial expansion.
But, since the following series are absolutely convergent and bounded,
$$\sum_{j=0}^\infty \left(\begin{array}{c} \alpha\\ j \end{array} \right) (-1)^{\alpha -j}h^{-itj}, ~~~~~~ \sum_{j=0}^\infty \left|\left(\begin{array}{c} \alpha\\ j \end{array} \right)\right|,$$
using the Abel- Stolz theorem for power series (see e.g. \cite{AH}), we obtain
$$\sum_{j=0}^\infty \left(\begin{array}{c} \alpha\\ j \end{array} \right) (-1)^{\alpha -j}h^{-itj} =
(h^{-it} - 1)^\alpha.$$
In order to justify the integration by series, we have for $s = c +it,$
\begin{eqnarray*}
&&\int_0^\infty |x^{s-1}| \sum_{j=0}^\infty \bigg|\left(\begin{array}{c} \alpha\\ j \end{array} \right)\bigg|h^{cj}|f(h^jx)|dx \\
&&= \sum_{j=0}^\infty \bigg|\left(\begin{array}{c} \alpha\\ j \end{array} \right)\bigg|h^{cj}|h^{j(1-s) - j}|\int_0^\infty |t^{s-1}| |f(t)| dt = \sum_{j=0}^\infty \bigg|\left(\begin{array}{c} \alpha\\ j \end{array} \right)\bigg| \| f\|_{X_c} < + \infty.
\end{eqnarray*}
As to iii), using i) and ii), taking the Mellin transform of both sides of the formula, we obtain
\begin{eqnarray*}
[\Delta_h^{\alpha,c}\Delta_h^{\beta,c}f]^\wedge_M(c+it) &=& (h^{-it} -1)^\alpha [\Delta_h^{\beta,c}f]^\wedge_M(c+it) =
(h^{-it} -1)^{\alpha  + \beta}[f]^\wedge_M(c+it)\\ &=& [\Delta_h^{\alpha + \beta,c}f]^\wedge_M(c+it),
\end{eqnarray*}
and so the assertion follows from the uniqueness theorem for Mellin transforms (see Theorem 8 in \cite{BJ}).
\vskip0,4cm
Note that the fractional differences introduced here depend fundamentally on the Mellin translation operator.
In the classical theories of Riemann-Liouville and Gr\"{u}nwald-Letnikov fractional calculus, the corresponding differences were based on the classical translation operator, and were first studied in a precise and 
systematic form in \cite{BW}; see also \cite{KST}, where  property iii) for these differences is also given, without proof. Moreover, other generlizations of fractional differences, via the Stirling functions of first kind, were also introduced in \cite{BKT4}, \cite{BKRT}.

\vskip0,4cm
\noindent
For spaces $X_{[a,b]},$ we have the following
\begin{Proposition}\label{Proposition 2}
Let $f \in X_{[a,b]},$ and let $c \in ]a,b[.$
\begin{description}
\item[(i)] If $0<h\leq 1,$ we have $\Delta_h^{\alpha,c}f \in X_{[a,c]},$ and for every $\nu \in [a,c[$
$$\|\Delta_h^{\alpha,c}f\|_{X_\nu} \leq \|f\|_{X_\nu} \sum_{j=0}^\infty \left|\left(\begin{array}{c} \alpha\\ j \end{array} \right)\right| h^{(c-\nu)j}.$$
Moreover,
$$M[\Delta_h^{\alpha,c}f](\nu +it) = (h^{c-\nu -it}-1)^\alpha M[f](\nu+it),~~t \in \mathbb{R}.$$
\item[(ii)] If $h> 1,$ we have $\Delta_h^{\alpha,c}f \in X_{[c,b]},$ and for every $\mu \in ]c,b]$
$$\|\Delta_h^{\alpha,c}f\|_{X_\mu} \leq \|f\|_{X_\mu} \sum_{j=0}^\infty \left|\left(\begin{array}{c} \alpha\\ j \end{array} \right)\right| h^{(c-\mu)j}.$$
Moreover,
\begin{eqnarray}
M[\Delta_h^{\alpha,c}f](\mu +it) = (h^{c-\mu -it}-1)^\alpha M[f](\mu+it),~~t \in \mathbb{R}.
\end{eqnarray}
\end{description}
\end{Proposition}
{\bf Proof}. We prove only (i) since the proof of (ii) is similar. Let $\nu \in [a,c[$ be fixed. Using an analogous reasoning as in Proposition \ref{Proposition 1}, we have
\begin{eqnarray*}
&&\|\Delta_h^{\alpha,c} f\|_{X_\nu}
\leq \int_0^{\infty} x^{\nu-1}\sum_{j=0}^\infty \bigg|\left(\begin{array}{c} \alpha\\ j \end{array} \right)\bigg|  h^{cj} |f(h^j x)|  dx
=  \|f\|_{X_\nu} \sum_{j=0}^\infty \bigg|\left(\begin{array}{c} \alpha\\ j \end{array} \right)\bigg|h^{j(c-\nu)},
\end{eqnarray*}
the last series being absolutely convergent for $0<h\leq 1.$
Moreover, as above, we can obtain, for $\nu \in [a,c[$ the assertion (5).~ $\Box.$
\vskip0,3cm
{\bf Definition 1}.
If for $f\in X_c$ there exists a function $g\in X_c$ such that
$$ \lim_{h\rightarrow 1}  \bigg\| \frac{\Delta_h^{\alpha,c} f(x)}{(h-1)^\alpha} - g(x) \bigg\|_{X_c} =0$$
then $g$ is called the strong fractional Mellin derivative of $f$ of order $\alpha,$ and it is denoted by $g(x)$ = s-$\Theta^\alpha_c f(x).$
If $\alpha =0$ it is easy to see that s-$\Theta^0_c f(x) =f(x).$
\vskip0,4cm
\noindent
We introduce now the Mellin Sobolev space $W^\alpha_{X_c}$ by
$$W^\alpha_{X_c} := \{ f\in X_c : \mbox{s-}\Theta^\alpha_c f~ \mbox{exists and}~  \mbox{s-}\Theta_c^\alpha f \in X_c \},$$
with $W^0_{X_c} = X_c.$
Analogously, for any interval $J$ we define the spaces $W^\alpha_{X_{J}},$ by
$$W^\alpha_{X_{J}} =\{f \in X_J: \mbox{s-}\Theta^\alpha_c f~ \mbox{exists for every}~ c \in J~\mbox{and}~\mbox{s-}\Theta^\alpha_c f\in X_J\}.$$
\vskip0,3cm
\noindent
For integral values of $\alpha$ our definition of strong Mellin derivative and the corresponding Mellin-Sobolev spaces reproduce those introduced in \cite{BJ1}, being the differences given now by a finite sum.

Using an approach introduced in \cite{BJ1} for the integer order case, we prove now
\newtheorem{Theorem}{Theorem}
\begin{Theorem}\label{Theorem 1}
The following properties hold:
\begin{description}
\item[(i)]
If $f\in W^\alpha_{X_c},$ then for $s=c+it, t\in \mathbb{R}$ we have
$$M[\mbox{s-}\Theta^\alpha_c f](s) = (-it)^\alpha M[f](s).$$
\item[(ii)] If $f\in W^\alpha_{X_{[a,b]}},$ then for every $\nu, c \in [a,b]$ we have
$$M[\mbox{s-}\Theta^\alpha_c f] (\nu +it) = (c-\nu-it)^\alpha M[f] (\nu + it)~~(t \in \mathbb{R}).$$
\end{description}
\end{Theorem}
{\bf Proof}.
As to (i), since
$$ \lim_{h\rightarrow 1} \bigg( \frac{h^{-it} -1}{h-1} \bigg)^\alpha = (-it)^\alpha,$$
we have, by Proposition \ref{Proposition 1}(ii),
\begin{eqnarray*}
\bigg | (-it)^\alpha [f]^{\wedge}_M (s) - [\mbox{s-}\Theta_c^\alpha f]^{\wedge}_M (s) \bigg| &=& \lim_{h\rightarrow 1} \bigg| \bigg( \frac{h^{-it} -1}{h-1} \bigg)^\alpha  [f]^{\wedge}_M (s) - [\mbox{s-}\Theta_c^\alpha f]^{\wedge}_M (s)\bigg|\\
&=& \lim_{h\rightarrow 1}\bigg|\bigg[\frac{\Delta_h^{\alpha,c} f}{(h-1)^\alpha}\bigg]^{\wedge}_M (s) - [\mbox{s-}\Theta_c^\alpha f]^{\wedge}_M (s)\bigg|\\
&=&
\lim_{h\rightarrow 1} \bigg| \bigg[\frac{\Delta_h^{\alpha,c} f}{(h-1)^\alpha} - \mbox{s-}\Theta_c^\alpha f \bigg]^{\wedge}_M(s) \bigg|\\ &\leq&
 \lim_{h\rightarrow 1} \bigg\| \frac{\Delta_h^{\alpha,c} f}{(h-1)^\alpha} - \mbox{s-}\Theta_c^\alpha f \bigg \|_{X_c} =0
\end{eqnarray*}
and thus (i) holds.
As to (ii), we can use the same approach, applying one-sided limits and Proposition \ref{Proposition 2}. $\Box$
\vskip0,4cm

\section{Mellin fractional integrals and the pointwise fractional Mellin differential operator}

In terms of Mellin analysis the natural operator of fractional integration is not the classical Liouville fractional integral of order $\alpha \in \mathbb{C},$ on $\mathbb{R}^+,$ with Re $\alpha >0,$ namely (1),
but the integral (2)
$$(J^\alpha_{0+} f)(x) = \frac{1}{\Gamma(\alpha)} \int_0^x  \bigg( \log \frac{x}{u} \bigg)^{\alpha-1} f(u) \frac{du}{u}~~~(x >0).$$
The above integrals were treated already in Mamedov's book \cite{MA}, page 168, in which the fractional integral of order $\alpha$ is defined by $(-1)^\alpha (J^\alpha_{0+} f)(x).$ This is due to a different notion of Mellin derivatives (of integral order), see Section 4.2. Our approach here is more direct and simple since it avoids the use of the complex coefficient $(-1)^\alpha.$

However,  for the development of the theory, it is important to consider the generalization of the  fractional integral, for $\mu \in \mathbb{R},$ in the form (3).

Note that for integer values $\alpha = r,$ in case $\mu = c$ and $f \in X_c$
(see \cite{BJ}, Definition 13), this turns into the iterated representation
\begin{equation}
(J^r_{0+,c}f)(x) = x^{-c}\int_0^x\int_0^{u_1}\ldots \int_0^{u_{r -1}} f(u_r) u_r^c \frac{du_r}{u_r}\ldots\frac{du_2}{u_2}\frac{du_1}{u_1}~~~(x >0).
\end{equation}

Several important properties of the operators $J^\alpha_{0+, \mu}$ were given by Butzer et al. in \cite{BKT}, \cite{BKT1},  \cite{BKT2}, (see also the recent monographs \cite{KST} and \cite{BDST}).
In particular, a boundedness property is given in the space $X_c,$  when the coefficient $\mu$ is greater than $c,$ (indeed a more general result is given there, for spaces $X^p_c$).
This is due to the fact that only for $\mu > c$ (or, in the complex case, Re $\mu > c$) we can view $J^\alpha_{0+, \mu} f$ as a Mellin convolution between two functions $f,g^\ast_\mu \in X_c,$ where
$$g^\ast_\mu(\frac{x}{u}):= (\frac{x}{u})^{-\mu}\frac{\chi_{]0,x]}(u)}{\Gamma (\alpha)}\bigg(\log (\frac{x}{u})\bigg)^{\alpha -1}.$$
Indeed, we have
\begin{eqnarray*}
&&(J_{0+,\mu}f)(x) = \frac{1}{\Gamma(\alpha)} \int_0^x  \bigg(\frac{u}{x}\bigg)^\mu\bigg(\log \frac{x}{u} \bigg)^{\alpha-1} f(u) \frac{du}{u} \\&=&
\frac{1}{\Gamma(\alpha)} \int_0^{+\infty}  \bigg(\frac{u}{x}\bigg)^\mu\chi_{]0,x]}(u) \bigg(\log \frac{x}{u} \bigg)^{\alpha-1} f(u) \frac{du}{u}\\
&=&\int_0^{+\infty} g^\ast_\mu(\frac{x}{u})f(u)\frac{du}{u} = (f\ast g^\ast_\mu)(x).
\end{eqnarray*}
Now, for $\mu >c$ the function:
$$g^\ast_\mu (u) = u^{-\mu}\frac{\chi_{]1,+\infty]}(u)}{\Gamma (\alpha)}(\log u)^{\alpha -1}$$
belongs to the space $X_c,$ as it is immediate to verify.

However, we are interested in properties of $J^\alpha_{0+, \mu} f$, when $\mu = c,$ since in the definition of the pointwise fractional Mellin derivative (see subsection 4.2), we have to compute such an integral with parameter $c.$ Hence in subsection 4.1 we will describe properties concerning the domain and the range of these fractional operators. As an example, we will show that for any non-trivial function $f$ in the domain of $J^\alpha_{0+, c}$ the image $J^\alpha_{0+, c}f$ cannot be in $X_c.$
This depends also on the fact that $g^\ast_c \not \in X_c.$
This implies that we cannot compute its Mellin transform of $g^\ast_c$ in the space $X_c$.

\subsection{The domain of $J^\alpha_{0+, c}$ and the semigroup property}

From now on we can consider the case $\alpha >0,$ the extension to complex $\alpha$ with Re $\alpha >0$ being similar but more technical.
 We define the domain of $J^\alpha_{0+, c},$ for $\alpha >0$ and $c \in \mathbb{R},$ as the class of all the functions $f:\mathbb{R}^+ \rightarrow \mathbb{C}$ such that
 \begin{eqnarray}
\int_0^x u^c \bigg( \log \frac{x}{u} \bigg)^{\alpha-1} |f(u)| \frac{du}{u} < + \infty
\end{eqnarray}
 for a.e. $x \in \mathbb{R}^+.$
 In the following we will denote the domain of $J^\alpha_{0+, c}$ by $Dom J^\alpha_{0+, c}$.

Recall that $X_{c, loc}$ is the space of all the functions such that $(\cdot)^{c-1}f(\cdot) \in L^1(]0,a[)$ for every $a >0.$
\begin{Proposition}\label{Proposition 3}
We have the following properties:
\begin{itemize}
\item[(i)] If $f \in X_{c,loc},$ then the function $(\cdot)^c f(\cdot) \in X_{1,loc}.$
\item[(ii)] If $c < c',$ then $X_{c,loc} \subset X_{c', loc}.$
\end{itemize}
\end{Proposition}
{\bf Proof.} (i) Let $a >0$ be fixed and let $f \in X_{c,loc}.$ Then
$$\int_0^a x^c |f(x)| dx = \int_0^a x x^{c-1}|f(x)| dx \leq a \int_0^a x^{c-1}|f(x)| dx$$
and so the assertion.

(ii) Let $f \in X_{c,loc}.$ Then, as before, setting $\alpha = c' -c,$ we can write
$$\int_0^a x^{c'-1}|f(x)| dx \leq a^{\alpha} \int_0^a x^{c-1}|f(x)|dx,$$
that is (ii) holds.
$\Box$
\vskip0,5cm
Note that the inclusion in (ii) does not hold for spaces $X_c.$
\vskip0,5cm
Concerning the domain of the operator $J^\alpha_{0+,c},$ we begin with the following proposition.
\begin{Proposition}\label{Proposition 4}
Let $\alpha >1,$ $c \in \mathbb{R}$ be fixed. Then $Dom J^\alpha_{0+,c} \subset X_{c,loc}.$
\end{Proposition}
{\bf Proof}. Assume that for a.e. $x \in \mathbb{R}^+$ the integral $(J^\alpha_{0+, c} |f|)(x),$ exists and put $F(u) = u^{c-1}f(u).$
We have to show that $F$ is integrable over $]0,a[,$ for any $a >0.$
Let $a >0$ be fixed and let $x > a$ be such that  $(J^\alpha_{0+, c} |f|)(x)$ exists. Then, for $u \in ]0,a[$ we have, since $|\log(x/u)|\leq |\log(x/a)|$
\begin{eqnarray*}
|F(u)| \leq \bigg| \log \frac{x}{u} \bigg|^{\alpha-1} |F(u)|  \bigg| \log \frac{x}{a} \bigg|^{1-\alpha},
\end{eqnarray*}
and the right-hand side of the inequality is integrable as a function of $u.$
$\Box$
\vskip0,4cm
Note that for $\alpha = 1$ we have immediately $DomJ^1_{0+,c} = X_{c,loc}.$

The case $0 < \alpha < 1$ is more delicate. We will show that in this instance $X_{c,loc} \subset Dom J^\alpha_{0+,c}.$

In order to give a more precise description of the domain of $J^\alpha_{0+, c},$ we now give a direct proof of the semigroup property in the domain of  fractional integrals.
This property is treated in \cite{BKT1}, \cite{KI} and \cite{KST}, but for the spaces $X^p_c(a,b)$ of all the functions $f: (a,b) \rightarrow \mathbb{C}$ such that
$(\cdot)^c f(\cdot) \in L^p(a,b),$ with $0< a < b \leq +\infty,\quad 1\leq p\leq \infty.$
However we prove this property under minimal assumptions, working directly in $Dom J^\alpha_{0+,c}$.

\begin{Theorem}\label{Theorem 2}
Let $\alpha, \beta >0, ~c \in \mathbb{R}$ be fixed. Let $f \in Dom J^{\alpha + \beta}_{0+,c}.$ Then
\begin{enumerate}
\item[(i)] $f \in Dom J^{\alpha}_{0+,c}\cap Dom J^{\beta}_{0+,c}$
\item[(ii)] $J^{\alpha}_{0+,c}f \in Dom J^{\beta}_{0+,c}$ and $J^{\beta}_{0+,c}f \in Dom J^{\alpha}_{0+,c}.$
\item[(iii)] $(J^{\alpha + \beta}_{0+,c}f)(x) = (J^{\alpha}_{0+,c} (J^{\beta}_{0+,c}f))(x),$~~a.e. $x \in \mathbb{R}^+.$
\item[(iv)] If $\alpha < \beta$ then $Dom J^\beta_{0+,c} \subset Dom J^\alpha_{0+,c}.$
\end{enumerate}
\end{Theorem}
{\bf Proof}. At first, let $f \in Dom J^{\alpha + \beta}_{0+,c}$ be a positive function. Then the integral
$$(J^{\alpha + \beta}_{0+, c}f)(x) = \frac{1}{\Gamma(\alpha + \beta)} \int_0^x \bigg( \frac{v}{x} \bigg)^c \bigg( \log \frac{x}{v} \bigg)^{\alpha + \beta -1} f(v) \frac{dv}{v}$$
is finite and nonnegative for a.e. $x \in \mathbb{R}^+.$

By  Tonelli's theorem on  iterated integrals of non-negative functions, and using  formula (2.8) concerning the Beta function in \cite{BKT1}, namely
$$\int_v^x  \bigg( \log  \frac{x}{u}\bigg)^{\alpha-1}\bigg( \log  \frac{u}{v}\bigg)^{\beta-1}\frac{du}{u} =
B(\beta, \alpha) \bigg( \log \frac{x}{v} \bigg)^{\alpha + \beta-1},$$
we have
\begin{eqnarray*}
&&(J^{\alpha + \beta}_{0+, c}f)(x) =
\frac{1}{\Gamma(\alpha) \Gamma (\beta)} \frac{\Gamma(\beta) \Gamma(\alpha)}{\Gamma(\alpha + \beta)} \int_0^x \bigg( \frac{v}{x} \bigg)^c \bigg( \log \frac{x}{v} \bigg)^{\alpha + \beta-1} f(v) \frac{dv}{v}\\ &=&
\frac{x^{-c}}{\Gamma(\alpha) \Gamma (\beta)}  \int_0^x v^c f(v) \bigg[ B(\beta, \alpha) \bigg( \log \frac{x}{v} \bigg)^{\alpha + \beta-1} \bigg] \frac{dv}{v}\\ &=&
\frac{x^{-c}}{\Gamma(\alpha) \Gamma (\beta)}  \int_0^{x} v^c  f(v) \bigg [\int_v^x  \bigg( \log  \frac{x}{u}\bigg)^{\alpha-1}\bigg( \log  \frac{u}{v}\bigg)^{\beta-1}\frac{du}{u} \bigg] \frac{dv}{v}\\ &=&
\frac{x^{-c}}{\Gamma(\alpha) \Gamma (\beta)} \int_0^x \int_0^{x} v^c \chi_{]v,x[}(u) \bigg( \log  \frac{x}{u}\bigg)^{\alpha-1}\bigg( \log  \frac{u}{v}\bigg)^{\beta-1} f(v)  \frac{dv}{v} \frac{du}{u}\\ &=&
\frac{x^{-c}}{\Gamma(\alpha)\Gamma (\beta)} \int_0^x \int_0^{x} v^c \chi_{]0,u[}(v) \bigg( \log  \frac{x}{u}\bigg)^{\alpha-1}\bigg( \log  \frac{u}{v}\bigg)^{\beta-1} f(v)  \frac{dv}{v} \frac{du}{u}\\ &=&
\frac{1}{\Gamma(\alpha)} \int_0^x   \bigg(   \frac{u}{x}\bigg)^c \bigg( \log  \frac{x}{u} \bigg)^{\alpha -1} \bigg[ \frac{1}{\Gamma(\beta)}\int_0^u \bigg(  \frac{v}{u} \bigg)^c \bigg(\log  \frac{u}{v}\bigg)^{\beta-1}  \frac{f(v)}{v} dv \bigg]  \frac{du}{u}\\&=& (J^{\alpha}_{0+, c}(J^\beta_{0+, c}f ))(x).
\end{eqnarray*}
This proves all the assertions (i), (ii), (iii), for positive functions.
In the general case, we can apply the above argument to the functions $f^+,~f^-$ using the linearity property  of the integrals.
Property (iv) follows immediately  by writing $\beta = \alpha + (\beta -\alpha)$ and applying (i).
$\Box$
\newtheorem{Corollary}{Corollary}
\begin{Corollary}\label{Corollary 1}
Let $0 < \alpha \leq 1, ~c \in \mathbb{R}$ be fixed. Then $X_{c,loc} \subset Dom J^\alpha_{0+,c}.$
\end{Corollary}
By this corollary, a consequence of (iv), we have the inclusions for
$\alpha < 1 < \beta,$
$$Dom J^\beta_{0+,c} \subset X_{c,loc} \subset Dom J^\alpha_{0+,c}.$$
These inclusions are strict. Indeed
\vskip0,3cm
\noindent
{\bf Examples}. For any $c \in \mathbb{R},~\beta >1,$ consider the function
$$f(x) = \frac{x^{-c}}{|\log x|^{\beta}}\chi_{]0,1/2[}(x).$$
Then $f \in X_{c,loc}$ but for any $x >1,$
\begin{eqnarray*}
&&\Gamma (\beta)(J^\beta_{0+,c}f)(x)= x^{-c}\int_0^x u^{c}\bigg(\log \frac{x}{u}\bigg)^{\beta -1}f(u) \frac{du}{u}=x^{-c}\int_0^{1/2}\frac{\bigg(\log \frac{x}{u}\bigg)^{\beta -1}}{u |\log u|^\beta}du\\
&\geq& x^{-c}\int_0^{1/2}\frac{\bigg(\log \frac{1}{u}\bigg)^{\beta -1}}{u |\log u|^\beta}du = x^{-c}\int_0^{1/2} \frac{1}{u |\log u|}du = + \infty.
\end{eqnarray*}
Moreover, for $0< \alpha < 1,$  consider the function:
\begin{eqnarray}
f(x) = \frac{x^{-c}}{|\log x|^{\gamma}}\chi_{]0,1/2[}(x),
\end{eqnarray}
where $\alpha <\gamma < 1.$ Then $f \not \in X_{c,loc},$ but for any $x > 1/2,$ we have
\begin{eqnarray*}
&&\Gamma (\alpha)(J^\alpha_{0+,c}f)(x) =
 x^{-c}\int_0^{1/2}\frac{1}{u \bigg(\log \frac{1}{u}\bigg)^{\gamma -\alpha +1}} \frac{\bigg(\log \frac{1}{u}\bigg)^{ 1-\alpha}}{\bigg(\log \frac{x}{u}\bigg)^{1-\alpha}}du \\
&\leq& \frac{M}{x^c}\int_0^{1/2} \frac{1}{u |\log u|^{\gamma -\alpha +1}}du < + \infty.
\end{eqnarray*}
\vskip0,4cm
Note that, more generally, the inclusion in (iv) of Theorem \ref{Theorem 2} is strict for any choice of $\alpha$ and $\beta.$
It is sufficent to consider the function  (8)
with $\alpha < \gamma < \beta.$
The calculations are the same.
\vskip0,4cm

We now give some sufficient conditions in order that a function $f$ belongs to the domain of the fractional integrals of order $\alpha >1.$
In this respect we have the following:

\begin{Proposition}\label{Proposition 5}
Let $\alpha >1.$
If $f\in X_{c,loc}$ is such that $f(u) = \mathcal{O}(u^{-(r +c -1)})$ for $u \rightarrow 0^+$ and $0<r<1,$ then
$f \in Dom J^\alpha_{0+,c}.$
\end{Proposition}
{\bf Proof.}  Let $x >0$ be fixed. Then we can write
$$ \int_0^x u^{c-1} |f(u)| \bigg(\log \frac{x}{u} \bigg)^{\alpha-1} du = \bigg( \int_0^{x/2} + \int_{x/2}^x \bigg) u^{c-1} |f(u)| \bigg(\log \frac{x}{u} \bigg)^{\alpha-1} du:= I_1 + I_2.$$
The integral $I_1$ can be estimated by considering the order of infinity at the point $0.$
The estimate of $I_2$ is easy since the function $\bigg(\log \frac{x}{u} \bigg)^{\alpha-1}$ is now bounded in the interval $[x/2,x].\Box$
\vskip0,4cm
Let us define
$$\widetilde{X}_{c,loc} =\{f \in X_{c,loc}: \exists r \in ]0,1[,~ \mbox{such that}~ f(u) = \mathcal{O}(u^{-(r +c-1)}),~u \rightarrow 0^+\}.$$
We have the following
\begin{Corollary}\label{Corollary 2}
Let $\alpha > 0,~c \in \mathbb{R}$ be fixed. Then
$$\widetilde{X}_{c,loc} \subset \bigcap_{\alpha >0} Dom J^\alpha_{0+,c}.$$
\end{Corollary}
\vskip0,4cm
Now let $f$ be a convergent power series of type
$$f(x) = \sum_{k=0}^\infty a_k x^k~~~~(a_k \in \mathbb{C},~k \in \mathbb{N}_0),$$
for $x \in [0,\ell],$ $\ell >0.$ For these functions the following series representation for $J^\alpha_{0+,c}f$ holds, when $c >0$
(see Lemma 4 and Lemma 5(i) in \cite{BKT3}):
$$(J^\alpha_{0+,c}f)(x) = \sum_{k=0}^{\infty}(c+k)^{-\alpha}a_kx^k~~~~(x \in [0,\ell]).$$
The assumption $c >0$ is essential. For $c=0,$ corresponding to the classical Hadamard integrals, we have the following
\begin{Proposition}\label{analitica}
Let $\alpha >0$ be fixed and let $f$ be a convergent power series as above. Then $f \in Dom J^\alpha_{0+}$ if and only if $f(0)=0.$ In this case we have
\begin{eqnarray}
(J^\alpha_{0+}f)(x) = \sum_{k=1}^\infty a_k k^{-\alpha}x^{k}~~~(0<x<\ell).
\end{eqnarray}
\end{Proposition}
{\bf Proof}. Let $f \in Dom J^\alpha_{0+}.$ Then the integral (7) is finite and
$$\int_0^x \bigg(\log \frac{x}{u}\bigg)^{\alpha -1}f(u)\frac{du}{u} = \int_0^x \bigg(\log \frac{x}{u}\bigg)^{\alpha -1}\sum_{k=1}^\infty a_k u^k\frac{du}{u} + a_0 \int_0^x \bigg(\log \frac{x}{u}\bigg)^{\alpha -1}\frac{du}{u} = I_1 +I_2.$$
As to $I_1$ we obtain
$$\int_0^x \bigg(\log \frac{x}{u}\bigg)^{\alpha -1}\sum_{k=1}^\infty |a_k| u^{k-1}du \leq \sum_{k=1}^\infty |a_k|x^{k-1}\int_0^x
 \bigg(\log \frac{x}{u}\bigg)^{\alpha -1}du.$$
 Since, using the change of variables $\log (x/u) = t,$
 $$\int_0^x \bigg(\log \frac{x}{u}\bigg)^{\alpha -1}du = x \Gamma(\alpha),$$
 we can integrate by series, yielding
 $$I_1 = \sum_{k=1}^\infty a_k \int_0^x \bigg(\log \frac{x}{u}\bigg)^{\alpha -1} u^{k-1}du < + \infty.$$
 As to $I_2,$ we get $I_2 < + \infty$ if and only if $a_0 = f(0) = 0,$ since
 $$\int_0^x\bigg(\log \frac{x}{u}\bigg)^{\alpha -1}\frac{du}{u} = + \infty.$$
 As to formula (9),
 $$(J^\alpha_{0+}f)(x) = \sum_{k=1}^\infty a_k \frac{1}{\Gamma(\alpha)}\int_0^x \bigg(\log \frac{x}{u}\bigg)^{\alpha -1} u^{k-1}du
 = \sum_{k=1}^\infty a_k (J^\alpha_{0+}t^k) (x) = \sum_{k=1}^\infty a_k k^{-\alpha}x^k,$$
 where in the last step we have applied the simple Lemma 3 in \cite{BKT3}, namely that $(J^\alpha_{0+}t^k)(x) = k^{-\alpha}x^k,$ $k>0.$ $\Box$
\vskip0,4cm
\noindent
Concerning the range of the operators $J^\alpha_{0+,c},$ we have the following important propositions.
\begin{Proposition}\label{Proposition 6}
Let $\alpha >0,~ c \in \mathbb{R}$ be fixed. If $f \in Dom J^{\alpha +1}_{0+,c},$ then $J^\alpha_{0+,c}f \in X_{c, loc}.$
\end{Proposition}
{\bf Proof}.
Let $f\in Dom J^{\alpha +1}_{0+,c}.$ We can assume that $f$ is nonnegative;
thus, for any $a >0,$
\begin{eqnarray*}
&&\Gamma (\alpha)\int_0^a x^{c-1}(J^{\alpha}_{0+,c}f)(x) dx
=\int_0^a u^{c-1} f(u)\bigg[\int_u^a \frac{1}{x}\bigg(\log \frac{x}{u}\bigg)^{\alpha -1}dx\bigg]du\\
&=&\frac{1}{\alpha} \int_0^a u^{c-1}f(u) \bigg(\log \frac{a}{u}\bigg)^{\alpha}du < + \infty.~~\Box
\end{eqnarray*}

\vskip0,3cm
Note that, in view of Proposition \ref{Proposition 6}, we can deduce that if $f \in Dom J^\alpha_{0+,c},$  not necessarily does
$J^\alpha_{0+,c}f \in X_{c,loc},$ unless $f \in Dom J^{\alpha +1}_{0+,c},$ which is a proper subspace of $Dom J^\alpha_{0+,c}.$
\vskip0,3cm
For example, we can take again the function $f$ of (8)
with $\alpha < \gamma < \alpha +1.$ Then $f \in Dom J^\alpha_{0+,c}$ but $f \not \in Dom J^{\alpha+1}_{0+,c}$ and
$J^\alpha_{0+, c}f \not \in X_{c,loc}.$
\vskip0,4cm

For spaces $X_c$ we have the following
\begin{Proposition}\label{Proposition 7}
Let $\alpha >0,~ c \in \mathbb{R}$ be fixed.
If $f\in Dom J^\alpha_{0+,c}$ is a non-negative function, then $J^\alpha_{0+, c} f \not \in X_{c},$ unless $f=0$ a.e. in $\mathbb{R}^+.$
\end{Proposition}
{\bf Proof}. Using an analogous argument as above, assuming $f \geq 0,$ we write
\begin{eqnarray*}
&&\int_0^{+\infty} x^{c-1} (J^{\alpha}_{0+, c}f)(x) dx
=\int_0^{+\infty} x^{-1} \bigg( \frac{1}{\Gamma(\alpha)} \int_0^{+\infty}  u^{c-1} \chi_{]0,x[}(u) \bigg( \log \frac{x}{u} \bigg)^{\alpha-1} f(u) du\bigg) dx\\ &=&
\int_0^{+\infty} \frac{1}{\Gamma(\alpha)} \bigg( \int_0^{+\infty}x^{-1} u^{c-1}  \chi_{]u, +\infty[}(x) \bigg( \log \frac{x}{u} \bigg)^{\alpha-1} f(u) dx\bigg) du\\ &=&
 \frac{1}{\Gamma(\alpha)} \int_0^{+\infty} \bigg( \int_u^{+\infty} x^{-1} \bigg( \log \frac{x}{u} \bigg)^{\alpha-1} dx \bigg) u^{c-1} f(u) du.
\end{eqnarray*}
Thus $(J^{\alpha}_{0+, c}f) \not \in X_c,$
since for every $u$
$$\int_u^{+\infty} \frac{1}{x ( \log \frac{x}{u})^{1-\alpha}} dx = +\infty. ~\Box$$
 \vskip0,3cm
 \noindent
The above result implies that a function $f \in Dom J^\alpha_{0+,c}$ such that $J^\alpha_{0+,c}f \in X_c,$ must necessarily change its sign. However the converse is not true in general, as proved by the following example: for a given $a >1,$ put
$$f(x) = -\chi_{[1/a, 1]}(x) + \chi_{]1, a]}(x).$$
It is easy to see that $f \in Dom J^\alpha_{0+,c} \cap X_c,$ but $J^\alpha_{0+,c}f \not \in X_c,$ for any $c \in \mathbb{R}.$
\vskip0,3cm
The following result, which will be useful in the following, is well known (see also \cite{BKT}, \cite{KST})
\begin{Proposition} \label{transform}
Let $\alpha >0$ and $c \in \mathbb{R}$ be fixed and let $f \in Dom J^\alpha_{0+,c}\cap X_c$ be such that $J^\alpha_{0+,c}f \in X_c.$ Then
$$M[J^\alpha_{0+,c}f](c +it) = (-it)^{-\alpha} M[f](c +it),~~t \in \mathbb{R}.$$
\end{Proposition}
\vskip0,4cm
Using Proposition \ref{transform} we study the structure of the functions $f$ such that $f \in Dom J^\alpha_{0+,c} \cap X_c$ for which $J^\alpha_{0+,c}f \in X_c.$
\begin{Proposition} \label{structure}
Let $f \in Dom J^\alpha_{0+,c} \cap X_c.$ If $J^\alpha_{0+,c}f \in X_c$ then
$$\int_0^{+\infty} x^{c-1}f(x) dx = 0.$$
\end{Proposition}
{\bf Proof}. Since $J^\alpha_{0+,c}f \in X_c,$ we can apply the Mellin transform on the line $s = c +it,$ obtaining
$$[J^\alpha_{0+,c}f]^\wedge_M(s) = (-it)^{-\alpha}[f]^\wedge_M(s) ~~(s= c+it, ~t \in \mathbb{R})$$
and this transform is a continuous and bounded function of $s.$ Therefore, taking $t=0$ we must have $[f]^\wedge_M(c) = 0,$ i.e. the assertion $\Box.$
\vskip0,4cm
Classes of functions $f \in Dom J^\alpha_{0+,c} \cap X_c$ for which $J^\alpha_{0+,c}f \in X_c$ may be easily constructed among the functions (of non-constant sign)
with compact support in $\mathbb{R}^+.$
\vskip0,4cm

However we have the following property (see also \cite{KST}, Lemma 2.33). We give the proof for the sake of completeness
\begin{Proposition}\label{Proposition 8}
Let $\alpha >0,~ c, \nu \in \mathbb{R},$  $\nu < c,$ being fixed.
If $f\in Dom J^\alpha_{0+,c}\cap X_{[\nu,c]},$ then $J^\alpha_{0+, c} f  \in X_{\nu}$ and
$$\|J^\alpha_{0+,c}f\|_{X_\nu} \leq \frac{\|f\|_{X_\nu}}{(c-\nu)^\alpha}.$$
Moreover, for any $s = \nu + it,$  we have
\begin{eqnarray*}
M[J^\alpha_{0+,c}f](\nu +it) = (c - \nu -it)^{-\alpha}M[f](\nu +it),~~t \in \mathbb{R}.
\end{eqnarray*}
$$|M[J^\alpha_{0+,c}f](s)| \leq \frac{\|f\|_{X_\nu}}{(c-\nu)^\alpha}.$$
\end{Proposition}
{\bf Proof}. We have by Tonelli's theorem,
\begin{eqnarray*}
&&\Gamma (\alpha) \|J^\alpha_{0+,c} f\|_{X_\nu} \\&\leq& \Gamma (\alpha) \int_0^{+ \infty} u^{\nu -1}|(J^\alpha_{0+,c}f)(u)| du \\&\leq&
\int_0^{+ \infty} u^{\nu -1}\bigg[\int_0^u \bigg(\frac{y}{u}\bigg)^c \bigg(\log \frac{u}{y}\bigg)^{\alpha -1} |f(y)| \frac{dy}{y}\bigg] du\\
&=&\int_0^{+\infty}\bigg[\int_y^{+\infty} u^{\nu -1 -c}\bigg(\log \frac{u}{y}\bigg)^{\alpha -1}du\bigg]y^{c-1}|f(y)|dy.
\end{eqnarray*}
For the inner integral, putting $\log (u/y) = z,$  we have:
$$\int_y^{+\infty} u^{\nu -1 -c}\bigg(\log \frac{u}{y}\bigg)^{\alpha -1}du =\int_0^{+\infty} y^{\nu-c} e^{-(c-\nu) z} z^{\alpha -1} dz =
\frac{y^{\nu -c}}{(c-\nu)^\alpha}\Gamma (\alpha),$$
and thus
$$\Gamma (\alpha)\|J^\alpha_{0+,c}f\|_{X_\nu} =
\frac{\Gamma (\alpha)}{(c-\nu)^\alpha}\|f\|_{X_\nu}.$$
As to the last part, the formula for the Mellin transform is established in \cite{BKT}, noting that
the Mellin transform on the line $s= \nu +it$ of the function $g^\ast_c(u) = u^{-c}(\log u)^{\alpha -1}\chi_{]1,+\infty[}(u) (\Gamma (\alpha))^{-1}$
is given by $[g^\ast_c](s) = (c-s)^{-\alpha} = (c-\nu -it)^{-\alpha},$
while for the estimate we easily have
$$|M[J^\alpha_{0+,c}f](s)| \leq \|J^\alpha_{0+,c}f\|_{X_\nu} = \frac{\|f\|_{X_\nu}}{(c-\nu)^\alpha}. \Box$$

\vskip0,4cm
Note that when $0<\alpha < 1,$ the assumption $f \in Dom J^\alpha_{0+,c}\cap X_{[\nu,c]},$ can be replaced by $f \in X_{[\nu,c]},$ since $X_{[\nu,c]} \subset Dom J^\alpha_{0+,c},$ by Corollary \ref{Corollary 1}.

\subsection{The pointwise fractional Mellin differential operator}

The pointwise fractional Mellin derivative of order $\alpha >0,$ or the Hadamard-type fractional derivative, associated with the integral $J^\alpha_{0+,c} f$,
$c \in \mathbb{R},$ and $f \in Dom J^{m- \alpha}_{0+, c},$ is given by
\begin{eqnarray}
(D^\alpha_{0+,c} f)(x) = x^{-c} \delta^m x^c (J^{m- \alpha}_{0+, c} f)(x)
\end{eqnarray}
where $m= [\alpha]+1$ and $\delta = \displaystyle(x \frac{d}{dx}).$ For $c=0,$ corresponding to the Hadamard fractional derivative, we put $(D^\alpha_{0+} f)(x):=(D^\alpha_{0+,0} f)(x).$
The above definition was introduced in \cite{BKT}, and then further developed in \cite{KI}, in which some sufficient conditions for the existence of the pointwise derivative are given in spaces of absolutely continuous type functions on bounded domains.
This notion originates from the theory of the classical Mellin differential operator, studied in \cite{BJ}. We give a short survey concerning this classical operator.

\vskip0,3cm
In the frame of Mellin transforms, the natural concept of a {\em pointwise} derivative of a function $f$ is given, as seen, by the limit of the difference quotient involving the Mellin translation; thus if $f'$ exists,
$$\lim_{h \rightarrow 1}\frac{\tau_h^cf(x) - f(x)}{h-1} = \lim_{h \rightarrow 1}\bigg[h^c x \frac{f(hx) - f(x)}{hx -x} + \frac{h^c -1}{h-1}f(x)\bigg] = x f'(x) + cf(x).$$
This gives the motivation of the following definition:
the pointwise Mellin differential operator $\Theta_c,$ or the pointwise Mellin derivative $\Theta_cf$ of a function $f: \mathbb{R}^+ \rightarrow \mathbb{C}$ and $c \in \mathbb{R},$ is defined by
\begin{eqnarray}
\Theta_cf(x) := x f'(x) + c f(x),~~x \in \mathbb{R}^+
\end{eqnarray}
provided $f'$ exists a.e. on $\mathbb{R}^+.$ The Mellin differential operator of order $r \in \mathbb{N}$ is defined iteratively by
\begin{eqnarray}
\Theta^1_c := \Theta_c ,\quad\quad \Theta^r_c := \Theta_c (\Theta_c^{r-1}).
\end{eqnarray}
For convenience set $\Theta^r:= \Theta^r_0$ for $c=0$ and $\Theta_c^0 := I,$ $I$ denoting the identity.
For instance, the first three Mellin derivatives are given by:
$$\Theta_cf(x) = xf'(x) + cf(x),$$
$$\Theta^2_cf(x) = x^2 f''(x) + (2c+1) xf'(x) + c^2f(x),$$
$$\Theta^3_cf(x) = x^3 f'''(x) + (3c+3)x^2f''(x) + (3c^2 + 3c +1)xf'(x) + c^3 f(x).$$

Let us return to Mamedov's book \cite{MA}. He defined the Mellin derivative of integral order 
in case $c=0$, in a slightly different, but essentially equivalent form, using the quotients
$$\frac{f(xh^{-1}) - f(x)}{\log h},$$
a definition connected with his notion of log-continuity.
It must be emphasised that this was a fully innovative  procedure at the time he introduced 
it. (His translation operator is actually $\tau_hf(x) = f(xh^{-1}),$ with incremental ratio $\log h,$ instead of $\log (1/h) = -\log 
h).$ His first order derivative, $(E^1f)(x)$, 
turns into, noting L'Hospital's rule,
$$(E^1 f)(x)= (-x)f'(x) =: -\Theta_0f(x).$$
  His derivatives of higher orders are defined inductively,
$$(E^mf)(x) = (-1)^m \Theta_0^mf(x),\quad m \in \mathbb{N}.$$
This is also Mamedov' s motivation of his definition of the fractional Mellin 
integral. Indeed, he used it to define the fractional derivative ($c=0$) for
  $ \alpha \in ]0,1[,$ by
$$(E^\alpha f)(x) := \lim_{h \rightarrow 1} 
\frac{(-1)^{1-\alpha} 
(J^{1-\alpha}_{0+}f)(xh^{-1}) - (J^{1-\alpha}_{0+}f)(x)]}{\log h},$$
and for $\alpha > 1,$ by
$$(E^\alpha f)(x) = E^{[\alpha]} (E^{\alpha -[\alpha]}f)(x).$$
Thus for example, if $\alpha \in ]0,1[,$, we have easily
$$(E^\alpha f)(x)  = 
(-1)^{2-\alpha}\Theta_0(J^{1-\alpha}_{0+}f)(x) = 
(-1)^{2-\alpha} (D^\alpha_{0+}f)(x),$$
which also gives  the link between Mamedov's 
definition and our present one. Analogously he proceeds  in case $\alpha >1.$
Using his definition of the fractional integral, 
Mamedov then studies the Mellin tranforms of the fractional integrals and derivatives of a function $f,$
(see section 23 of \cite{MA}). 
From these results it would have been possible to deduce a version of the fundamental theorem of the integral and 
differential calculus in his fractional frame, in 
the special case when the function $f,$ its fractional derivative and  fractional integral 
belong to the space $X_0.$ However he presents it explicitely only for integer values of $\alpha,$ 
(formula (22.3), page 169).
Nevertheless it is indeed a surprising result, the only comparative result being that for  the Chebyshev transform
  \cite{BS} of 1975. For this very reason is the  late Prof. Mamedov a true pioneer of Mellin analysis.
  On the other hand, the approach given in \cite{BJ} is somewhat more direct and simpler, and  the 
 present  versions of the fundamental theorem in  local spaces $X_{c,loc},$ , given in Theorems 
 \ref{Theorem 3} and \ref{Theorem 3bis} below,  are more general and elegant. 
But recall that Mamedov's first papers appeared in 1979/81,   \cite{MO0, MO1, MO2},  thus almost twenty years
 earlier than \cite{BJ}.
\vskip0,4cm
We have the following
\begin{Proposition}\label{Proposition 9}
We have, for $m \in \mathbb{N},~x>0,$
$$\delta^m x^c f(x) = x^c \Theta_c^m f (x).$$
\end{Proposition}
{\bf Proof.}
For $m=1$ we have
$$\delta x^c  f(x) = x(c x^{c-1} f(x) + x^c f' (x)) = x^c (c f(x) + x f'(x)) = x^c \Theta_c f(x).$$
Now we suppose that the relation holds for $m$ and prove that it holds for $m+1.$
$$ \delta^{m+1}(x^c f(x)) = \delta (\delta^m (x^c f(x)) = \delta (x^c \Theta_c^m f(x)) = x^c \Theta_c(\Theta^m_c f(x))= x^c  \Theta_c^{m+1}f(x),$$
and so the assertion.$\Box$
\vskip0,4cm
For $r \in \mathbb{N},$ $\Theta^r_cf(x)$ is given by the following proposition, also giving the connections between Mellin and ordinary derivatives
(these relations was also given in \cite{BKT}, \cite{KST}, but without proofs).
\begin{Proposition}\label{Proposition 10}
Let $f \in X_{c,loc}$ be such that $\Theta^r_cf(x)$ exists at the point $x$ for $r \in \mathbb{N}.$ Then $(D^r_{0+,c} f)(x)$ exists and
$$(D^r_{0+,c} f)(x) = \Theta^r_cf(x) = \sum_{k=0}^r S_c(r,k) x^k f^{(k)}(x),$$
where $S_c(r,k),$ $0\leq k\leq r,$ denote the generalized Stirling numbers of second kind, defined recursively by
$$S_c(r,0) := c^r,~ S_c(r,r) := 1,~ S_c(r+1,k) = S_c(r, k-1)+ (c+k) S_c(r,k).$$
 In particular for $c=0$
 $$\Theta^r f(x) = \sum_{k=0}^r S(r,k) x^k f^{(k)}(x)$$
 $S(r,k):= S_0(r,k)$ being the (classical) Stirling numbers of the second kind.
\end{Proposition}
{\bf Proof}. For $r = 1$ (that is $m=2$), we have
\begin{eqnarray*}
&&(D^1_{0+,c}f)(x)
=x^{-c}\delta^2 x^c \frac{1}{\Gamma (1)} \int_0^x \bigg(\frac{u}{x}\bigg)^c \bigg(\log \frac{x}{u}\bigg)^{1-1} f(u) \frac{du}{u}\\
&&= x^{-c}\delta \bigg(x \frac{d}{dx}\bigg) \int_0^x u^{c-1}f(u) du = x^{-c}\delta x^c f(x) = \Theta_cf(x).
\end{eqnarray*}
\noindent
For $r=2,$ (that is $m=3$), we have
\begin{eqnarray*}
&&(D^2_{0+,c}f)(x)
=x^{-c}\delta^3 x^c \frac{1}{\Gamma (2)} \int_0^x \bigg(\frac{u}{x}\bigg)^c \bigg(\log \frac{x}{u}\bigg)^{1-1} f(u) \frac{du}{u}\\
&&= x^{-c}\delta \bigg(x \frac{d}{dx}\bigg)x^c f(x) =x^{-c}\delta x(cx^{c-1} f(x) + x^c f'(x))\\
&&= cxf'(x) + c^2f(x) + (c+1)xf'(x) + x^2 f''(x) = x^2 f''(x) + (2c+1) xf'(x) + c^2 f(x)\\&&=
\Theta^2_cf(x).
\end{eqnarray*}
In the general case, using Proposition \ref{Proposition 9} we have
\begin{eqnarray*}
&&(D^{r}_{0+,c}f)(x)
= x^{-c}\delta^{r}(\delta x^c (J^1_{0+,c}f)(x))
\\&&= x^{-c} \delta^{r}(x^c f(x)) =  x^{-c} x^c \Theta^r_cf(x) = \Theta^r_cf(x).
\end{eqnarray*}
Now in accordance with (11) and (12) we have (see \cite{BJ})
\begin{eqnarray*}
&&\Theta^{r+1}_cf(x) = \Theta_c (\Theta^r_cf)(x) = x\frac{d}{dx}\Theta^r_cf(x) + c \Theta^r_cf(x)\\
&&=\sum_{k=0}^r S_c(r,k)((k+c)x^k f^{(k)}(x) + x^{k+1}f^{(k+1)}f(x)) =
\sum_{k=0}^{r+1}S_c(r+1,k)x^k f^{(k)}(x),
\end{eqnarray*}
and so the assertion follows. $\Box$
\vskip0,4cm
\noindent
Note that in Proposition \ref{Proposition 10} the basic assumption that $f \in X_{c,loc}$ is essential. Let for example $g(x) = 1,$ for every
$x \in \mathbb{R}^+$ and $c=0.$ Then $g \not \in X_{0,loc} = Dom J^1_{0+}.$ This implies that we cannot compute $D^r_{0+}f,$ while obviously we have
$\Theta^rf(x) = 0,$ for any $r \in \mathbb{N}.$
Another example is given by the function $h(x) = \log x,$ $x \in \mathbb{R}^+.$ In this instance, for $c=0$ and $r=1$ we have $\Theta h(x) = 1,$
while $h \not \in X_{0, loc}.$
\vskip0,4cm
\noindent
Now we turn to the fractional case. The above Proposition shows that the notion of Hadamard-type fractional derivative $D^\alpha_{0+,c}$ is the natural extension of the Mellin derivative $\Theta_c^kf,$ with $k \in \mathbb{N},$ to the fractional case as also applies to the ordinary and Riemann-Liouville fractional derivatives.
A simple consequence of Proposition \ref{Proposition 9} is the following alternative representation of the fractional derivative of $f$, for $\alpha >0$
$$(D^\alpha_{0+,c} f)(x) = \Theta_c^{m} (J^{m- \alpha}_{0+, c} f)(x)$$
where $m= [\alpha]+1.$
Using this representation we can obtain the following Proposition
\begin{Proposition}\label{Prop. carlo}
Let $\alpha>0, c \in \mathbb{R},$ be fixed and $m-1 \leq \alpha < m.$ Let $f\in X_{c, loc}$ be such that $f^{(m)}\in X_{c, loc},$ then
$$ (D^\alpha_{0+,c} f)(x) =\sum_{k=0}^m S_c(m,k) x^k (J^{m- \alpha}_{0+, c+k} f^{(k)})(x)).$$
\end{Proposition}
{\bf Proof.}
At first note that from the assumptions, for any $0<\gamma \leq 1$ the derivatives $f^{(k)},$ $k=1,\ldots m,$ belongs to the domain of
$J^\gamma_{0+,c+k}.$
Note that using a simple change of variable we can write, for every $c \in \mathbb{R},$
$$(J^\gamma_{0+,c}f)(x) = \frac{1}{\Gamma(\gamma)}\int_1^{+ \infty}\frac{1}{v^{c+1}}(\log v)^{\gamma -1} f(\frac{x}{v})dv.$$
Thus differentiating under the integral we easily have
$$(J^\gamma_{0+,c}f)'(x) = \frac{1}{\Gamma(\gamma)}\int_1^{+ \infty}\frac{1}{v^{c+2}}(\log v)^{\gamma -1} f'(\frac{x}{v})dv =
(J^\gamma_{0+,c+1}f')(x)$$
and by an easy induction we obtain, for $x >0$ and $k \in \mathbb{N},$
$$(J^\gamma_{0+,c}f)^{(k)}(x) = (J^\gamma_{0+,c+k}f^{(k)})(x).$$
Hence by Lemma 9 in \cite{BJ}, (see also Proposition \ref{Proposition 10}), we have
$$(D^\alpha_{0+,c}f)(x) = \Theta^m_c(J^{m-\alpha}_{0+,c}f)(x) =
\sum_{k=0}^m S_c(m,k) x^k (J^{m-\alpha}_{0+, c+k}f^{(k)})(x),$$
that is the assertion $\Box.$
\vskip0,4cm

First, let $f$ be a convergent power series as in Proposition \ref{analitica}. In this instance, we obtain the following formula
for the derivative $D^\alpha_{0+}f$:
\begin{Proposition}\label{analitica 2}
Let $\alpha >0$ be fixed and $f$ be as in Proposition \ref{analitica}, such that $f(0) =0.$ Then for $0<x<\ell,$
$$(D^\alpha_{0+}f)(x) = \sum_{k=1}^\infty a_k k^\alpha x^k.$$
\end{Proposition}
{\bf Proof}. Putting $m = [\alpha] +1,$ by integration and differentiation by series, using similar reasonings as in Proposition \ref{analitica},
we have
\begin{eqnarray*}
&&(D^\alpha_{0+}f)(x) = \delta^m \frac{1}{\Gamma (m-\alpha)}\int_0^x \bigg(\log \frac{x}{u}\bigg)^{m-\alpha -1}
\sum_{k=1}^\infty a_k u^k \frac{du}{u}= \delta^m \sum_{k=1}^\infty a_k (J^{m -\alpha}_{0+}t^k)(x)\\
 &=& \delta^m \sum_{k=1}^\infty a_k k^{-(m-\alpha)}x^k
= \sum_{k=1}^\infty a_k k^{-(m-\alpha)}\delta^m x^k = \sum_{k=1}^\infty a_k k^{\alpha}x^k.~\Box
\end{eqnarray*}
\vskip0,3cm
The above Proposition extends Lemma 5 (ii) in \cite{BKT3} to the case $c=0.$
\vskip0,3cm
An interesting representation, for analytic functions, of the derivative $D^\alpha_{0+,c}f$ is given in terms
of infinite series involving the  Stirling functions of the second kind $S_c(\alpha, k),$ which can be defined for $c \in \mathbb{R}$ by
$$S_c(\alpha, k) := \frac{1}{k!}\sum_{j=0}^k (-1)^{k-j} \left(\begin{array}{c} k\\ j \end{array} \right)(c+j)^\alpha~~~(\alpha \in \mathbb{C},~k \in I\!\!N_0).$$
This representation, given in \cite{BKT3}, is as follows:
\begin{Proposition}\label{analitica 3}
Let $f:\mathbb{R}^+ \rightarrow \mathbb{R}$ be an arbitrarily often differentiable function such that its Taylor series converges and let
$\alpha > 0, ~c >0.$ Then
$$(D^\alpha_{0+,c}f)(x) = \sum_{k=0}^\infty S_c(\alpha, k)x^k f^{(k)}(x)~~~(x>0).$$
\end{Proposition}
For $c=0$ also an inverse formula is available, expressing the classical Riemann-Liouville fractional derivative in terms of the Mellin derivatives (see \cite{BKRT}), namely
$$x^\alpha({\mathcal D}^\alpha_{0+}f)(x) = \sum_{k=0}^\infty s(\alpha, k)(D^k_{0+}f)(x),~~(\alpha >0,~x>0),$$
where ${\mathcal D}^\alpha_{0+}f$ denotes the Riemann-Liouville fractional derivative and $s(\alpha, k)$ the Stirling functions of the first kind.

An analogous representation holds also for the fractional integrals $J^\alpha_{0+,c}f,$ namely (see \cite{BKT3})
\begin{Proposition}\label{analitica 4}
Let $f \in Dom J^\alpha_{0+,c},$ and $f:\mathbb{R}^+ \rightarrow \mathbb{R}$ satisfy the hyphothesis of Proposition 13. Then
$$(J^\alpha_{0+,c}f)(x) = \sum_{k=0}^\infty S_c(-\alpha, k)x^k f^{(k)}(x)~~~(x>0).$$
\end{Proposition}
\vskip0,4cm
Since $(D^\alpha_{0+,c}f)(x),~(J^\alpha_{0+,c}f)(x),$ for $\alpha >0,$ and $S_c(\alpha, k),$ for $\alpha \in \mathbb{R},~k \in \mathbb{N}_0,$ are three continuous functions of $c \in \mathbb{R}$ at $c=0,$ we can let $c\rightarrow 0$ in the previous Propositions, and can deduce corresponding representations of Hadamard fractional differentiation and integrations in terms of the Stirling functions $S(\alpha,k)$ and classical derivatives, if both $J^\alpha_{0+}f$ and $D^\alpha_{0+}f$ exists (for details see \cite{BKT3}).
\vskip0,4cm
\noindent
Now we introduce certain Mellin-Sobolev type spaces which will be useful in the following (see also \cite{BJ}). Firstly,
we define
$$AC_{loc}:=\{f: \mathbb{R}^+\rightarrow \mathbb{C}: f(x) = \int_0^xg(t)dt,~ \mbox{for a given} ~ g\in L^1_{loc}(\mathbb{R}^+)\}.$$

 Recall that $L^1_{loc}(\mathbb{R}^+)$ stands for the space of all (Lebesgue) measurable functions $g:\mathbb{R}^+ \rightarrow \mathbb{C}$ such that
 $$\int_0^x g(t)dt$$
 exists as a Lebesgue integral for every $x >0.$
 For any $f \in AC_{loc}$ we have $f' = g$ a.e., where $f'$ denotes the usual derivative.
For any $c \in \mathbb{R},$ we define
$$AC_{c,loc} :=\{f\in X_{c,loc}: (\cdot)^c f(\cdot) \in AC_{loc}\}.$$
   For any $c\in \mathbb{R}$  we define $AC^1_{c,loc} = AC_{c,loc}$ and for $m \in \mathbb{N},$ $m\geq 2,$
 $$AC_{c,loc}^m :=\{f \in AC_{c,loc}: \delta^{m-1}((\cdot)^cf(\cdot)) \in AC_{loc}\}.$$
 
 We have the following
 \begin{Lemma} \label{absolutecont}
 If $f \in AC_{c,loc}^m,$ then the Mellin derivative $\Theta_c^mf$ exists and $\Theta_c^mf \in X_{c,loc}.$
 \end{Lemma}
 {\bf Proof}. Since $\delta^{m-1}((\cdot)^cf(\cdot)) \in AC_{loc},$ we have
 $$\frac{d}{dx}\delta^{m-1}(x^cf(x)) \in L^1_{loc}(\mathbb{R}^+).$$
 But, using Proposition \ref{Proposition 9}
 $$\frac{d}{dx}\delta^{m-1}(x^cf(x)) = x^{-1}\delta^m(x^c f(x)) = x^{c-1}\Theta_c^mf(x),$$
 and so the assertion follows. $\Box$
 \vskip0,4cm
 \begin{Lemma} \label{asintotico}
 If $f \in AC^m_{c,loc},$ $m \geq 2,$ then $\delta^j((\cdot)^c f(\cdot)) \in AC_{loc},$ for $j=0,1,\ldots,m-2,$ and
  $$\lim_{x \rightarrow 0^+}\delta^j((x)^c f(x)) = 0.$$
 \end{Lemma}
 {\bf Proof}. The case $m=2$ follows immediately from the definitions, while for $m>2$ one can use the relation
 $$ \delta^{j-1}((x)^cf(x)) = \int_0^x \delta^j((u)^cf(u))\frac{du}{u},~~j=1,2\ldots m-2. \Box$$
 \vskip0,4cm
 The following result gives a representation of functions in $AC^m_{c,loc}.$ A similar result for functions defined on a compact interval $[a,b] \subset \mathbb{R}^+$ is given in \cite{KI}.
 \begin{Lemma} \label{representation}
 Let $f \in AC^m_{c,loc},$ $m\geq 1,$ and let us assume that $\varphi_m:= \frac{d}{dx}\delta^{m-1}((\cdot)^c f(\cdot))
 \in Dom J^m_{0+,1}.$ If there exists $\alpha \in ]0,1[$ such that  $\varphi_m (x) = {\cal O}(x^{-\alpha}),~~x \rightarrow 0^+,$
 then we have necessarily
 $$f(x) = x^{1-c}J^m_{0+,1}\varphi_m(x).$$
 \end{Lemma}
 {\bf Proof}. For $m=1,$ there exists $\varphi_1 \in L^1_{loc}(\mathbb{R}^+)$ such that
 $$f(x) = x^{-c}\int_0^x \varphi_1(t)dt = x^{1-c}J^1_{0+,1}\varphi_1(x),$$
 and so the assertion follows.

 For $m=2,$ $\delta((\cdot)^cf(\cdot)) \in AC_{loc},$ and so there exists $\varphi_2 \in L^1_{loc}(\mathbb{R}^+)$ such that
 $$\delta(t^cf(t)) = \int_0^t \varphi_2(u)du.$$
 Let $\varepsilon >0$ be fixed. Integrating the above relation in the interval $[\varepsilon, x]$ we have
 $$\int_\varepsilon^x \delta(t^cf(t))\frac{dt}{t} = \int_\varepsilon^x\bigg(\int_0^t \varphi_2(u)du\bigg)\frac{dt}{t}.$$
 Integrating by parts, we get
 \begin{eqnarray*}
 &&x^cf(x) - \varepsilon^c f(\varepsilon) = \bigg[\log t \int_0^t \varphi_2(u)du \bigg]_\varepsilon^x -
 \int_\varepsilon^x \log t \varphi_2(t) dt\\
 &=&\log x \int_0^\varepsilon \varphi_2(t) dt + \int_\varepsilon^x \log \frac{x}{t}~\varphi_2(t) dt -
 \log \varepsilon \int_0^\varepsilon \varphi_2(t) dt.
 \end{eqnarray*}
 Letting $\varepsilon \rightarrow 0^+,$ since $\varphi_2 \in Dom J^2_{0+,1},$ by Lemma 4, we obtain
 $$x^cf(x) = \int_0^x \log \frac{x}{t}~\varphi_2(t) dt - \lim_{\varepsilon \rightarrow 0^+} \log \varepsilon \int_0^\varepsilon \varphi_2(t) dt.$$
 Since by assumption, $\varphi_2(t) = {\cal O}(t^{-\alpha}),~~t \rightarrow 0^+,$ using the De L'Hopital rule, the limit on the right-hand side of the previous relation is zero. Thus,
 $$f(x) = x^{-c} \int_0^x \log \frac{x}{t}~\varphi_2(t) dt = x^{1-c}J^2_{0+,1}\varphi_2(x).$$
 For the general case one can apply the same method, using the binomial formula $\Box$

  \vskip0,4cm
 Now, for every $c \in \mathbb{R}$ and $m \in \mathbb{N},$ we introduce the Mellin-Sobolev space by
 $${\mathcal X}_{c,loc}^m :=\{f \in X_{c,loc}: f=g ~\mbox{a.e. in}~ \mathbb{R}^+, \mbox{for}~g \in AC^m_{c,loc}\}$$

  A non-local version of the above space, denoted by ${\mathcal X}_{c}^m$ is defined in \cite{BJ}. It contains
 all the functions $f:\mathbb{R}^+ \rightarrow \mathbb{C}$ such that $f \in X_c$ and there exists $g \in
 AC^m_{c,loc}$ such that $f=g$ a.e. in $\mathbb{R}^+$ with $\Theta^m_cf \in X_c.$

 Note that, if $f \in X_c$ is such that $J^m_{0+,c}f \in X_c$ then $J^m_{0+,c}f \in {\mathcal X}^m_c$ (see \cite{BJ}, Theorem 11).

  In particular, a function $f \in {\mathcal X}_{c}^m$ is such that $f\in X_c,$  $\Theta^m_{c}f$ exists and $\Theta^m_{c}f \in X_c.$
This suggests a way to define the fractional versions of the above spaces.
For a given $\alpha >0,$
we define
$${\mathcal X}_{c}^\alpha :=\{f \in X_{c}: (D^\alpha_{0+,c} f)(x)~\mbox{exists a.e. and}~ D^\alpha_{0+,c} f \in X_{c}\}$$
and its local version
$${\mathcal X}_{c,loc}^\alpha :=\{f \in X_{c,loc}: (D^\alpha_{0+,c} f)(x)~\mbox{exists a.e. and}~ D^\alpha_{0+,c} f \in X_{c,loc}\}.$$
Analogously we can define the spaces ${\mathcal X}^\alpha_{J},$ for any interval $J$ as
$${\mathcal X}^\alpha_{J}:=\{f \in X_J: (D^\alpha_{0+,c} f)(x)~ \mbox{exists a.e. for every}~ c \in J~\mbox{and}~(D^\alpha_{0+,c} f)(x)\in X_J\}$$
and its local version ${\mathcal X}^\alpha_{J, loc}.$
We begin with the following
\begin{Proposition}\label{Proposition 11}
 Let $f \in {\mathcal X}^\alpha_{c,loc}$ be such that $\Theta_c^mf \in X_{c,loc},$ where $m = [\alpha] + 1.$ Then
$$(D^\alpha_{0+,c} f)(x)= \Theta_c^m(J^{m-\alpha}_{0+,c}f)(x) = J^{m-\alpha}_{0+,c} (\Theta_c^m f) (x).$$
\end{Proposition}
{\bf Proof.}
Since $f\in X_{c,loc}$ and $0<m-\alpha <1,$ $f \in Dom J^{m-\alpha}_{0+,c}$ and $\Theta^m_cf \in Dom J^{m-\alpha}_{0+,c}$
by Corollary \ref{Corollary 1}.
The first equality is already stated as a consequence of  Proposition \ref{Proposition 9}, thus we will prove the other equality.
We obtain  by (10)
\begin{eqnarray*}
&&(D^\alpha_{0+,c}f)(x) =  x^{-c} \delta^m (x^c (J^{m-\alpha}_{0+,c}f))(x)\\
&&=x^{-c} \bigg(\delta^m \bigg[x^c \frac{1}{\Gamma (m-\alpha)}\int_0^x
\bigg(\frac{v}{x}\bigg)^c \bigg(\log\frac{x}{v}\bigg)^{m -\alpha-1} f(v) \frac{dv}{v}\bigg]\bigg)(x)\\
&&=  x^{-c}\bigg(\delta^m \bigg[ \frac{1}{\Gamma (m-\alpha)}\int_1^{+ \infty} \frac{x^c}{t^{c+1}}(\log t)^{m-\alpha-1}f(\frac{x}{t})dt\bigg]\bigg)(x) \\
&&=  x^{-c}\sum_{k=0}^m S(m,k) x^k \frac{d^k}{dx^k}  \bigg[\frac{1}{\Gamma (m-\alpha)}\int_1^{+ \infty}
\frac{x^c}{t^{c+1}}(\log t)^{m-\alpha-1}f(\frac{x}{t}) dt \bigg]\\
&&= \frac{x^{-c}}{\Gamma (m-\alpha)}\int_1^{+ \infty} \sum_{k=0}^m S(m,k) x^k \frac{d^k}{dx^k}(x^c f(\frac{x}{t})) (\log t)^{m-\alpha-1} \frac{dt}{t^{c+1}}.
\end{eqnarray*}
Using the elementary formula for the derivatives of the product, we have
\begin{eqnarray*}
&&(D^\alpha_{0+,c}f)(x) \\
&&= \frac{x^{-c}}{\Gamma (m-\alpha)}\int_1^{+ \infty} \sum_{k=0}^m S(m,k) x^k \sum_{j=0}^k \left(\begin{array}{c} k\\ j \end{array} \right)
\prod_{\nu =0}^{j-1}(c-\nu)\frac{x^{c-j}}{t^{k-j}}f^{(k-j)}(x/t)(\log t)^{m-\alpha-1} \frac{dt}{t^{c+1}}\\
&&=
\frac{x^{-c}}{\Gamma (m-\alpha)}\int_0^x \sum_{k=0}^m S(m,k) v^k \frac{d^k}{dv^k} (v^c f(v))(\log (x/v))^{m-\alpha-1} \frac{dv}{v}\\
&&= \frac{x^{-c}}{\Gamma (m-\alpha)}\int_0^x (\delta^m (v^cf(v))(\log (x/v))^{m-\alpha-1} \frac{dv}{v}= J^{m-\alpha}_{0+,c} (\Theta_c^m f) (x).
 \end{eqnarray*}
 where we have used Proposition \ref{Proposition 9}. Thus the assertion follows. $\Box$
 \vskip0,4cm
 In order to prove a new fractional version of the fundamental theorem of the differential and integral calculus in the Mellin frame, first we give the following proposition concerning the case $\alpha = m\in \mathbb{N}.$
 Recall that in this case, using the representation in terms of iterated integrals, $J^m_{0+,c}f$ is m-times differentiable, whenever $f \in Dom J^m_{0+,c}.$
 \begin{Proposition}\label{Proposition 12}
 We have:
 \begin{enumerate}
 \item[(i)] Let $f \in {\mathcal X}^1_{c,loc},$ then
 $$J^1_{0+,c}(\Theta_cf)(x) = f(x),~~a.e.~ x \in \mathbb{R}^+.$$
 \item[(ii)] Let $m \in \mathbb{N}, m>1,$ and let $f \in {\mathcal X}^m_{c, loc}$ be such that $\Theta^m_cf \in Dom J^m_{0+,c}.$ Then
 $$J^m_{0+,c}(\Theta^m_cf)(x) = f(x),~~a.e.~ x \in \mathbb{R}^+.$$
 \item[(iii)] Let $f \in X_{c,loc},$ then
 $$\Theta_c (J^1_{0+,c}f)(x) = f(x),~~a.e.~ x \in \mathbb{R}^+.$$
 \item[(iv)] Let $f \in Dom J^m_{0+,c},$ then
 $$\Theta^m_c (J^m_{0+,c}f)(x) = f(x),~~a.e. ~x \in \mathbb{R}^+.$$
 \end{enumerate}
 \end{Proposition}
 {\bf Proof}. As to (i) we have, by the absolute continuity and Lemma \ref{asintotico}
 $$J^1_{0+,c}(\Theta_cf)(x) = \int_0^x \bigg(\frac{u}{x}\bigg)^c (\Theta_cf) (u)\frac{du}{u} =
 x^{-c}\int_0^x \frac{d}{du}(u^c f(u))du = f(x),$$
 a.e. $x \in \mathbb{R}^+.$

 For (ii) we can obtain the result using the iterated representation of $J^m_{0+,c}f,$ m-times the absolute continuity, Lemma \ref{asintotico} and Proposition \ref{Proposition 9}.
 
 As to (iii) note
 $$\Theta_c(J^1_{0+,c}f)(x) = x \frac{d}{dx} (J^1_{0+,c}f)(x) + c (J^1_{0+,c}f)(x).$$
 and  we have
 $$x \frac{d}{dx} (J^1_{0+,c}f)(x) =
  -cx^{-c}\int_0^x u^{c-1}f(u)du + f(x)$$
 from which we obtain the assertion.

 Finally we prove (iv) and we will use again induction.
 Assuming that (iv) holds for $m-1,$ we have by Theorem \ref{Theorem 2} (iii),
 \begin{eqnarray*}
 \Theta^{m}_c(J^{m}_{0+,c}f)(x) &=& \Theta_c (\Theta^{m-1}_c (J^{m-1}_{0+,c}(J^{1}_{0+,c}f)))(x) = f(x),
 \end{eqnarray*}
 a.e. $x \in \mathbb{R}^+,$ by the induction assumption. $\Box$
 \vskip0,4cm
Proposition \ref{Proposition 12} gives a version of Theorem 11 in \cite{BJ} for the  spaces $X_{c,loc},$ without the use of Mellin transforms and under sharp assumptions. As a consequence, for spaces $X_c,$ we deduce again the formula for the Mellin transform
of $J^m_{0+,c}f$ whenever $J^m_{0+,c}f \in X_c$
$$[J^m_{0+,c}f]^\wedge_M(c+it) = (-it)^{-m}[f]^\wedge_M(c+it).$$
\vskip0,4cm
 Now we are ready to prove the fundamental theorem of the fractional differential and integral calculus in the Mellin frame.

\begin{Theorem}\label{Theorem 3}
Let $\alpha >0$ be fixed and $m= [\alpha] + 1.$
 \begin{enumerate}
 \item[a)]
 Let $f \in {\mathcal X}^\alpha_{c,loc} \cap {\mathcal X}^m_{c, loc},$ such that $D^\alpha_{0+,c}f, \Theta^{m}_cf \in Dom J^{m}_{0+,c}.$
 Then
 $$(J^\alpha_{0+, c}(D^\alpha_{0+,c}f))(x) = f(x),~a.e. ~x \in \mathbb{R}^+.$$
 \item[b)]
 Let $f \in Dom J^{m}_{0+,c},$ be  such that  $J^\alpha_{0+,c}f \in X_{c,loc}.$ Then
$$(D^\alpha_{0+,c}(J^\alpha_{0+,c}f))(x) = f(x),~~a.e. ~x \in \mathbb{R}^+.$$
\end{enumerate}
 \end{Theorem}
 {\bf Proof}. As to part a), by Propositions \ref{Proposition 9}, \ref{Proposition 11}, \ref{Proposition 12} and Theorem \ref{Theorem 2}, we have
 for a.e. $x \in \mathbb{R}^+$
\begin{eqnarray*}
&&(J^\alpha_{0+, c}(D^\alpha_{0+,c}f))(x) = J_{0+, c}^\alpha \bigg(x^{-c}(\delta^m(x^c J_{0+,c}^{m- \alpha}f))\bigg)(x) \\&&
=(J^\alpha_{0+, c}(J_{0+,c}^{m- \alpha}(\Theta_c^mf)))(x)
= (J^m_{0+, c}(\Theta_c^mf))(x) = f(x).
\end{eqnarray*}
As to part b), we have, by Propositions \ref{Proposition 9}, \ref{Proposition 12} and Theorem \ref{Theorem 2},
\begin{eqnarray*}
&&(D^\alpha_{0+,c}(J_{0+,c}^\alpha f))(x) = x^{-c}\delta^m(x^c J_{0+,c}^{m-\alpha}(J_{0+,c}^\alpha f))(x)\\
&&= x^{-c}\delta^m (x^c J_{0+,c}^m f)(x) = \Theta_c^m (J_{0+,c}^mf)(x) = f(x),
\end{eqnarray*}
almost everywhere. $\Box$
\vskip0,4cm
\noindent
A result related to part a) is also described in \cite{KST}, Lemma 2.35, for  functions $f$ belonging to the subspace
$$J^\alpha_{0+, \mu}(X^p_c):=\{f= J^\alpha_{0+,\mu}g,~\mbox{for}~ g \in X^p_c\}$$
with $\mu > c.$ In this instance the formula is a simple consequence of part b), using the integral representation of $f.$
Related results in spaces $X^p_\nu$ with $c> \nu$ are given in \cite{KST}, Property 2.28. Note that
for $p=1,$ if $f \in X_\nu,$ with $c >\nu,$ then $J^\alpha_{0+,c}f \in X_\nu,$ so that our assumption is satisfied.
For bounded intervals $I$ similar results are also given in \cite{KI}, for functions belonging to $L^p(I)$.
\vskip0,4cm
\noindent
More generally, we can, with our approach, also consider compositions between the operators of Hadamard-type fractional integrals and derivatives, in local spaces (for similar results in $X^p_c$ spaces see \cite{KI} on bounded intervals, and \cite{KST} in $\mathbb{R}^+$).
\begin{Theorem}\label{Theorem 3bis}
Let $\alpha, \beta >0$ with $\beta>\alpha$ and $m= [\alpha] + 1.$
\begin{enumerate}
\item[a')]
Let $f \in {\mathcal X}^\alpha_{c,loc}\cap {\mathcal X}^m_{c, loc},$ such that $D^\alpha_{0+,c}f \in Dom J^{\beta}_{0+,c}$ and
 $\Theta^{m}_cf \in Dom J^{m+\beta -\alpha}_{0+,c}.$
 Then
 $$(J^\beta_{0+, c}(D^\alpha_{0+,c}f))(x) = (J^{\beta -\alpha}_{0+,c}f)(x),~a.e. ~x \in \mathbb{R}^+.$$
 \item[b')]
 Let $f \in Dom J^{m + \beta - \alpha}_{0+,c}.$  Then
  $$(D^\alpha_{0+,c}(J^\beta_{0+, c}f))(x) = (J^{\beta - \alpha}_{0+, c}f)(x),~a.e. ~x \in \mathbb{R}^+.$$
 \end{enumerate}
 \end{Theorem}
 {\bf Proof}. Regarding a'), as in the proof of Theorem \ref{Theorem 3},
  we have  for a.e. $x \in \mathbb{R}^+$
\begin{eqnarray*}
&&(J^\beta_{0+, c}(D^\alpha_{0+,c}f))(x)
=(J^\beta_{0+, c}(J_{0+,c}^{m- \alpha}(\Theta_c^mf)))(x)
= J^{\beta -\alpha}_{0+,c}(J^m_{0+, c}(\Theta_c^mf))(x) = (J^{\beta -\alpha}_{0+,c}f)(x).
\end{eqnarray*}
Regarding b'), we have
 \begin{eqnarray*}
&&(D^\alpha_{0+,c}(J_{0+,c}^\beta f))(x)
= x^{-c}\delta^m (x^c J_{0+,c}^{m-\alpha +\beta} f)(x) = \Theta_c^m (J_{0+,c}^{m +\beta -\alpha} f)(x) =
(J^{\beta - \alpha}_{0+,c}f)(x).~~ \Box
\end{eqnarray*}
\vskip0,4cm

\section{A relation between pointwise and strong fractional Mellin derivatives}

In this section we will compare the definitions of the Mellin derivative in the strong and pointwise versions.
For this purpose, we need some further notations and preliminary results.

For $h,x \in \mathbb{R}^+$ and $\widetilde{c} \in \mathbb{R},$ we define
\begin{eqnarray*}
m^{\widetilde{c}}_h(x) : = \left\{\begin{array}{ll} x^{-\widetilde{c}}\chi_{[1/h,1]}(x),~~\mbox{if}~h\geq 1,\\
-x^{-\widetilde{c}}\chi_{[1,1/h]}(x),~~\mbox{if}~0<h <1.
\end{array} \right.
\end{eqnarray*}
It is clear that $m^{\widetilde{c}}_h \in X_{(-\infty, \infty)}$ and we have (see \cite{BJ})
\begin{eqnarray*}
[m^{\widetilde{c}}_h]^\wedge_M(s) = \left\{\begin{array}{ll} \frac{1}{\widetilde{c} -s} (h^{\widetilde{c} -s}-1),~~\mbox{if}~s\in \mathbb{C}\setminus \{\widetilde{c}\},\\
\log h,~~\mbox{if}~s = \widetilde{c}.
\end{array} \right.
\end{eqnarray*}
Denoting the $r$th-fold convolution of $m^{\widetilde{c}}_h$ with itself by $(m^{\widetilde{c}}_h \ast)^r,~r \in \mathbb{N},$ we have, by Theorem 3 in \cite{BJ}
$$[(m^{\widetilde{c}}_h \ast)^r]^\wedge_M(s) = (\widetilde{c}-s)^{-r}(h^{\widetilde{c}-s}-1)^r,~~s \in \mathbb{C}\setminus \{\widetilde{c}\}.$$
We recall that by Proposition \ref{Proposition 2}, one has for $f \in X_{[a,b]},$  $c,\nu \in [a,b]$ and $r \in \mathbb{N},$
\begin{eqnarray}
M\left[\sum_{k=0}^r (-1)^{r-k} \left(\begin{array}{c} r\\ k \end{array} \right)\tau^{c}_{h^k}f\right](\nu + it) =
(h^{c-\nu-it} -1)^r M[f](\nu + it).
\end{eqnarray}
\vskip0,3cm
In \cite{BJ} (Proposition 6, formula (8.8)), the following lemma was established:
\begin{Lemma}\label{Lemma 2}
If $f \in {\mathcal X}^r_{[a,b]},~r \in \mathbb{N},$ then for $c \in [a,b],~h>1$ we have, for $x \in \mathbb{R}^+$
$$\sum_{k=0}^r (-1)^{r-k} \left(\begin{array}{c} r\\ k \end{array} \right)\tau^{c}_{h^k}f(x) =
x^{-c}\bigg((m^0_h\ast)^r \ast (\Theta^r_cf(\cdot) (\cdot)^c)\bigg)(x).$$
\end{Lemma}
\vskip0,3cm
\begin{Lemma}\label{Lemma 3}
If $f \in {\mathcal X}^r_{[a,b]},~r \in \mathbb{N},$ then for $c, \nu \in [a,b],$ we have
$$M[\Theta^r_cf](\nu +it) = (c- \nu -it)^r M[f](\nu +it).$$
\end{Lemma}
{\bf Proof}. Let us put, for $x \in \mathbb{R}^+$
$$G(x) = \bigg((m^0_h\ast)^r \ast (\Theta^r_cf(\cdot) (\cdot)^c)\bigg)(x).$$
Since $\Theta^r_cf \in X_\nu$ by assumption, it is easy to see that $\Theta^r_cf(\cdot)(\cdot)^c \in X_{\nu -c}.$ Then
$G \in X_{\nu -c}$ and so $(\cdot)^{-c}G(\cdot) \in X_\nu.$
Hence by Lemma \ref{Lemma 2}, Proposition \ref{Proposition 2} and (13) we have
\begin{eqnarray*}
&&M[(\cdot)^{-c}G(\cdot)](\nu +it) = M\left[\sum_{k=0}^r (-1)^{r-k} \left(\begin{array}{c} r\\ k \end{array} \right)\tau^{c}_{h^k}f\right](\nu +it) \\
&&= (h^{c-\nu -it}-1)^r M[f](\nu +it).
\end{eqnarray*}
Using Proposition 1(c) in \cite{BJ}  and the convolution theorem, (Lemma \ref{convolution}(ii)),
we have
\begin{eqnarray*}
M[(\cdot)^{-c}G(\cdot)](\nu +it) &=& M[G](\nu -c +it)  \\
&=&M[(m^0_h\ast)^r](\nu -c +it)
M[\Theta^r_cf(\cdot)(\cdot)^c](\nu-c+it)\\
&=& (h^{c -\nu -it} -1)^r (c-\nu -it)^{-r}M[\Theta^r_cf](\nu +it),
\end{eqnarray*}
from which we deduce the assertion. $\Box$
\vskip0,4cm
 We prove the main theorem of this section.
\begin{Theorem}\label{Theorem 5}
Let $\alpha >0$ be fixed.
\begin{enumerate}
\item[(i)] Let $f \in {\mathcal X}^\alpha_{c}$ such that $\Theta^{m}_c f \in X_{c},$ where $m = [\alpha] +1.$ Then $f \in W^\alpha_{X_{c}}$ and
$$(D^\alpha_{0+,c}f) (x) = \mbox{s-}\Theta^\alpha_cf(x),~~a.e. ~x \in \mathbb{R}^+.$$
\item[(ii)] Let $f \in {\mathcal X}^\alpha_{[a,b]}$ such that $\Theta^{m}_c f \in X_{[a,b]},$ for $c \in ]a,b[,$ where $m = [\alpha] +1.$ Then $f \in W^\alpha_{X_{[a,b]}}$ and
$$(D^\alpha_{0+,c}f) (x) = \mbox{s-}\Theta^\alpha_cf(x),~~a.e. ~x \in \mathbb{R}^+,~~~c \in ]a,b[.$$
\end{enumerate}
\end{Theorem}
{\bf Proof}. (i) By Proposition \ref{Proposition 11} we have
$$(D^\alpha_{0+,c}f) (x) = (J^{m -\alpha}_{0+,c}(\Theta^{m}_cf))(x),$$
which belongs to $X_c.$
Thus, passing to Mellin transforms, we have, for $t \in \mathbb{R},$
\begin{eqnarray*}
&&[D^\alpha_{0+,c}f ]^\wedge_M(c + it) = [(J^{m -\alpha}_{0+,c}(\Theta^{m}_cf))]^\wedge_M(c +it) \\
&&= (-it)^{\alpha-m}[\Theta_c^{m}f]^\wedge_M(c +it) = (-it)^\alpha [f]^\wedge_M(c +it) = [\mbox{s-}\Theta^\alpha_cf]^\wedge_M(c+it).
\end{eqnarray*}
Hence, $D^\alpha_{0+,c}f$ and s-$\Theta^\alpha_cf$ have the same Mellin transform along the line $s = c +it$, and so
the assertion follows by the identity theorem (see \cite{BJ}). \\
(ii) Again, using Proposition \ref{Proposition 11}, and taking the Mellin transform on the line $s=\nu + it,$ for $\nu \in ]a,b[$ with $\nu < c,$ we obtain
$$[D^\alpha_{0+,c}f ]^\wedge_M(\nu + it) = [J^{m-\alpha}_{0+,c}(\Theta^m_cf)]^\wedge_M(\nu + it) = (c-\nu -it)^\alpha[f]^\wedge_M(\nu +it),$$
and so the assertion follows as before. $\Box$
\vskip0,3cm
The above theorem reproduces Theorem 4.3 in \cite{BJ1}, for integral values of $\alpha,$ i.e. for $k \in \mathbb{N},$ the pointwise Mellin derivative $\Theta^k_cf$ equals the strong derivative s-$\Theta^k_cf,$ as defined in \cite{BJ1}, for functions belonging to the space $W^k_{X_c}.$

As a consequence of Theorem \ref{Theorem 5}, for the spaces $X_c$ we can give more direct proofs of the fundamental formulae of integral and differential calculus in the fractional Mellin setting, now using the Mellin transform. We begin with the following
\begin{Theorem}\label{Theorem 3ter}
Let $\alpha >0$ be fixed.
\begin{enumerate}
\item[a)]
Let $f \in {\mathcal X}^\alpha_c$ be such that $D^\alpha_{0+,c}f \in Dom J^\alpha_{0+,c},$ and $J^\alpha_{0+,c}f \in X_c.$ Then
$$(J^\alpha_{0+, c}(D^\alpha_{0+,c}f))(x) = f(x)~a.e. ~x \in \mathbb{R}^+.$$
\item[b)]
Let $f \in Dom J^{\alpha}_{0+,c}\cap X_{c}$ be such that  $J^{\alpha}_{0+,c}f \in {\mathcal X}^\alpha_c.$ Then we have
$$(D^\alpha_{0+,c}J^\alpha_{0+, c}f)(x) = f(x),~a.e. ~x \in \mathbb{R}^+.$$
\end{enumerate}
\end{Theorem}
 {\bf Proof}. As to part a), we can compute the Mellin transforms, obtaining
 \begin{eqnarray*}
 &&[J^\alpha_{0+, c}(D^\alpha_{0+,c}f)]^\wedge_M(c +it) = (-it)^{-\alpha}[D^\alpha_{0+,c}f]^\wedge_M(c +it)\\
 &=&
 (-it)^{-\alpha}(-it)^{\alpha}[f]^\wedge_M(c +it) = [f]^\wedge_M(c +it)
 \end{eqnarray*}
 and so the assertion follows by the uniqueness theorem of Mellin transform.

 Part b) is carried out using the same approach. $\Box$
 \vskip0,5cm
 In comparison with Theorem \ref{Theorem 3bis} we have, under different assumptions, the following
\begin{Theorem}\label{Theorem 3quater}
Let $\alpha, \beta >0$ with $\beta > \alpha.$
\begin{enumerate}
\item[a')]
Let $f \in {\mathcal X}^\alpha_{c}.$  If
 $D^\alpha_{0+,c}f \in Dom J^\beta_{0+,c},$ and $J^\beta_{0+, c}(D^\alpha_{0+,c}f) \in X_c,$ then
 $$(J^\beta_{0+, c}(D^\alpha_{0+,c}f))(x) = J^{\beta -\alpha}_{0+,c}f(x),~a.e. ~x \in \mathbb{R}^+.$$
 \item[b')]
Let $f \in Dom J^{\beta}_{0+,c}\cap X_{c}$ be such that  $J^{\beta}_{0+,c}f \in {\mathcal X}^\alpha_c.$ Then
  $$(D^\alpha_{0+,c}J^\beta_{0+, c}f)(x) = (J^{\beta - \alpha}_{0+, c}f)(x),~a.e. ~x \in \mathbb{R}^+.$$
 \end{enumerate}
 \end{Theorem}
 {\bf Proof}. As to part a'), using again the Mellin transform, we have
  \begin{eqnarray*}
 &&[J^\beta_{0+, c}(D^\alpha_{0+,c}f)]^\wedge_M(c +it) = (-it)^{-\beta}[D^\alpha_{0+,c}]^\wedge_M(c +it)\\
 &=&
 (-it)^{-\beta}(-it)^{\alpha}[f]^\wedge_M(c +it) = [J^{\beta -\alpha}_{0+,c}f]^\wedge_M(c +it),
 \end{eqnarray*}
 and so the assertion follows again by the uniqueness theorem.
 As to part b'), the proof is similar to the previous one. $\Box$

\vskip0.3cm
\noindent
For the special case of spaces $X_{[a,b]},$ we have the following two further results.
  \begin{Theorem}\label{Theorem 3quint}
Let $\alpha >0$ be fixed.
\begin{enumerate}
\item[a'')]
 Let $f \in {\mathcal X}^\alpha_{[a,b]}$ and $c \in ]a,b].$ If
 $D^\alpha_{0+,c}f \in Dom J^\alpha_{0+,c},$ then
 $$(J^\alpha_{0+, c}(D^\alpha_{0+,c}f))(x) = f(x) ~a.e. ~x \in \mathbb{R}^+.$$
 \item[b'')]
 Let $c,\nu \in [a,b]$ with $\nu < c.$ If $f \in Dom J^{\alpha}_{0+,c}\cap X_{[a,b]}$ is such that  $J^{\alpha}_{0+,c}f \in {\mathcal X}^\alpha_\nu,$ then
  $$(D^\alpha_{0+,c}J^\alpha_{0+, c}f)(x) = f(x),~a.e. ~x \in \mathbb{R}^+.$$
 \end{enumerate}
 \end{Theorem}
 {\bf Proof}. As to part a'') take $\nu \in [a,b]$ with $\nu < c.$ Using the Mellin transform on the line $s = \nu +it,$ we have by Proposition \ref{Proposition 8} and Theorem \ref{Theorem 5}
 \begin{eqnarray*}
 &&[J^\alpha_{0+, c}(D^\alpha_{0+,c}f)]^\wedge_M(\nu +it) = (c-\nu -it)^{-\alpha}[D^\alpha_{0+,c}]^\wedge_M(\nu +it)\\
 &=&
 (c- \nu -it)^{-\alpha}(c - \nu -it)^{\alpha}[f]^\wedge_M(\nu +it) = [f]^\wedge_M(\nu +it)
 \end{eqnarray*}
 and so the assertion again follows by the uniqueness theorem.  Part b'') is carried out similarly. $\Box$
  \vskip0,4cm
 \begin{Theorem}\label{Theorem 3sest}
Let $\alpha, \beta >0$ be fixed with $\beta> \alpha.$
\begin{enumerate}
\item[a''')]
 Let $f \in {\mathcal X}^\alpha_{[a,b]}$ and $c \in ]a,b].$ If
 $D^\alpha_{0+,c}f \in Dom J^\beta_{0+,c},$ then
 $$(J^\beta_{0+, c}(D^\alpha_{0+,c}f))(x) = (J^{\beta - \alpha}_{0+,c}f)(x) ~a.e. ~x \in \mathbb{R}^+.$$
 \item[b''')]
 Let $c,\nu \in [a,b]$ with $\nu < c.$ If $f \in Dom J^{\beta}_{0+,c}\cap X_{[a,b]}$ is such that  $J^{\beta}_{0+,c}f \in {\mathcal X}^\alpha_\nu,$ then
  $$(D^\alpha_{0+,c}J^\beta_{0+, c}f)(x) = (J^{\beta - \alpha}_{0+,c}f)(x),~a.e. ~x \in \mathbb{R}^+.$$
 \end{enumerate}
 \end{Theorem}
 {\bf Proof}. The proof is essentially the same as in Theorem \ref{Theorem 3quint} taking Mellin transforms in the space $X_\nu.$ $\Box$
 \vskip0,4cm
 For what concerns the strong fractional Mellin derivatives we have
\begin{Theorem}\label{Theorem 7} Let $\alpha>0$ be fixed.
\begin{enumerate}
\item[(i)] Let $f \in W^\alpha_{X_c}$ be such that $\mbox{s-}\Theta^\alpha_cf \in Dom J^\alpha_{0+,c}$ and $J^\alpha_{0+,c}(\mbox{s-}\Theta^\alpha_cf) \in X_c.$ Then
$$J^\alpha_{0+,c}(\mbox{s-}\Theta^\alpha_cf)(x) = f(x),~~\mbox{a.e}~x \in \mathbb{R}^+.$$
\item[(ii)] Let $f \in W^\alpha_{X_{[a,b]}}$ be such that $\mbox{s-}\Theta^\alpha_cf \in Dom J^\alpha_{0+,c},$ for $c \in ]a,b[.$ Then
$$J^\alpha_{0+,c}(\mbox{s-}\Theta^\alpha_cf)(x) = f(x),~~\mbox{a.e}~x \in \mathbb{R}^+.$$
\end{enumerate}
\end{Theorem}
{\bf Proof}. (i) By assumptions, we can compute the Mellin transform of the function $J^\alpha_{0+,c}(\mbox{s-}\Theta^\alpha_cf),$ on the line $s=c+it,$ obtaining, as before,
$$[J^\alpha_{0+,c}(\mbox{s-}\Theta^\alpha_cf)]^\wedge_M (c +it) = [f]^\wedge_M(c +it),~~(t \in \mathbb{R}).$$
Analogously (ii) follows, by taking Mellin transforms on $s=\nu +it,$ with $\nu < c.$
$\Box$
\vskip0,3cm
\begin{Theorem}\label{Theorem 8}
Let $\alpha >0$ be fixed.
\begin{enumerate}
\item[(i)] Let $f \in X_c$ be such that $J^\alpha_{0+,c}f \in W^\alpha_{X_c}.$ Then
$$\mbox{s-}\Theta^\alpha_c(J^\alpha_{0+,c}f)(x) = f(x),~~\mbox{a.e}~x \in \mathbb{R}^+.$$
\item[(ii)] Let $f \in X_{[a,b]}$ be such that $J^\alpha_{0+,c}f \in W^\alpha_{X_{\nu}},$ $\alpha >0$ and $c \in ]a,b], \nu < c.$ Then
$$\mbox{s-}\Theta^\alpha_c(J^\alpha_{0+,c}f)(x) = f(x),~~\mbox{a.e}~x \in \mathbb{R}^+.$$
\end{enumerate}
\end{Theorem}
{\bf Proof}. The proof is essentially the same as for the previous theorem. $\Box$
\vskip0,4cm
\noindent
In order to state an extension to the fractional setting of the equivalence theorem proved in \cite{BJ} Theorem 10, we introduce the following
subspace of $W^\alpha_{X_c},$ for $\alpha >0$ and $c \in \mathbb{R},$
$$\widetilde{W}^\alpha_{X_c} =\{f \in W^\alpha_{X_c} : \mbox{s-}\Theta^\alpha_cf \in Dom J^\alpha_{0+,c} ~\mbox{and}~ J^\alpha_{0+,c}(\mbox{s-}\Theta^\alpha_cf) \in X_c\}.$$
\begin{Theorem}\label{Theorem 6}
Let $f\in X_{c}$ and $\alpha>0.$
The following four assertions are equivalent
\begin{itemize}
\item[(i)] $f\in \widetilde{W}^\alpha_{X_{c}}$.
\item[(ii)] There is a function $g_1\in  X_{c}\cap Dom J^\alpha_{0+,c}$ with $J^\alpha_{0+,c}g_1 \in X_c$ such that
$$ \lim_{h\rightarrow 1} \bigg\| \frac{\Delta^{\alpha,c}_h f}{(h-1)^\alpha} - g_1 \bigg \|_{X_c} =0.$$
\item[(iii)] There is $g_2\in  X_{c} \cap Dom J^\alpha_{0+,c}$ with $J^\alpha_{0+,c}g_2 \in X_c$ such that
$$ (-it)^\alpha M[f](c+it) = M[g_2](c+it).$$
\item[(iv)] There is $g_3\in X_{c}\cap Dom J^\alpha_{0+,c}$ such that $J^\alpha_{0+,c}g_3 \in X_c,$ and
$$ f(x) = \frac{1}{\Gamma(\alpha)} \int_0^x  \bigg(\frac{u}{x}\bigg)^c \bigg( \log \frac{u}{x}\bigg)^{\alpha-1} g_3(u) \frac{du}{u}\quad a.e. \quad x\in \mathbb{R^+}.$$
\end{itemize}
If one of the above assertions is satisfied, then $D^\alpha_{0+,c} f(x) =$ s-$\Theta^\alpha_c f(x)= g_1 = g_2= g_3$ a.e. $x\in \mathbb{R^+}.$
\end{Theorem}
{\bf Proof.}
It is easy to see that (i) implies (ii) and (ii) implies (iii) by Theorem \ref{Theorem 1}.
We prove now (iii) implies (iv).
Let $g_2\in X_{c}$ be such that (iii) holds.
Then, putting $g_3= g_2$
 we have, by Proposition \ref{transform},
$$M[J^\alpha_{0+,c} g_2](c+it) = (-it)^{-\alpha}M[g_3](c+it) = M[f](c + it).$$
Thus by (iii) we have immediately the assertion by the identity theorem for Mellin transforms.
Finally we prove that (iv) implies (i). By (iv), we have in particular that $J^\alpha_{0+,c}g_3 \in X_{c,loc}.$ This implies
that $J^\alpha_{0+,c}g_3 \in Dom J^{m-\alpha}_{0+,c},$ since $0 <m-\alpha < 1.$ So, by the semigroup property (Theorem \ref{Theorem 2}),
$g_3 \in Dom J^m_{0+,c}.$ Therefore the assumptions of Theorem \ref{Theorem 3} part b) are satisfied and we have
$$ (D^\alpha_{0+,c}f)(x) = (D^\alpha_{0+,c} (J^\alpha_{0+,c} g_3))(x) = g_3(x) \quad a.e. \quad x\in \mathbb{R^+}.$$
So the assertion follows. $\Box$
\vskip0.3cm
\noindent
Analogous equivalence theorems hold for the spaces $W^\alpha_{X_{[a,b]}}.$

\section{Some particular applications}

In this section we discuss some basic examples.
\begin{enumerate}
\item[1.] The first example also discussed in \cite{BKT3}, Lemma 3 and in Property 2.25 in \cite{KST} and used in the proof of Propositions \ref{analitica} and \ref{analitica 2}, is the following:
Consider the function $g(x) = x^b,$ $b \in \mathbb{R}.$ Then for any $c \in \mathbb{R}^+$ such that $c+b > 0$ we have $g \in Dom J^\alpha_{0+,c},$ and
$$(J^\alpha_{0+,c}g)(x) = (c+b)^{-\alpha}x^b.$$
 In particular, for $b>0$ and $c=0$ we get $(J^\alpha_{0+}g)(x) = b^{-\alpha}x^b.$ Analogously,
we have also
$$(D^\alpha_{0+,c}g)(x) = (c+b)^{\alpha}x^b.$$
 This also well enlightens the fundamental theorem in the fractional frame.

It should be noted that in this case we cannot compute $J^\alpha_{0+}1$ and $D^\alpha_{0+}1$ since the function $g(t) = 1,$ corresponding to $b=0,$
is not in the domain of $J^\alpha_{0+}.$ However we can compute $J^\alpha_{0+,c}1$ and $D^\alpha_{0+,c} 1,$ with $c>0,$ obtaining easily $(J^\alpha_{0+,c} 1)(x) = c^{-\alpha},$ and  $(D^\alpha_{0+,c} 1)(x) = c^\alpha.$ The last relation follows by
$\delta^m(x^c c^{\alpha - m}) = c^\alpha x^c,$ for $m-1 < \alpha < m,$ which is proved by an easy induction.

Moreover we could also calculate $J^\alpha_{a+}1$ and $D^\alpha_{a+}1,$ with $a>0$ in place of $0$ in the definitions of the Hadamard-type integrals and derivatives (see \cite{KST}).

\item[2.] As a second example, let us consider the function $g_k(x) = \log^kx,$ for $k \in \mathbb{N}.$ For any $\alpha >0$ and $c>0$ we have
$g_k \in DomJ^\alpha_{0+,c}$ and by a change of variables and using the binomial theorem, we can write
\begin{eqnarray*}
(J^\alpha_{0+,c}g_k)(x) &=& \frac{1}{\Gamma (\alpha)}\int_0^{+\infty} e^{-cv} v^{\alpha -1}(\log x - v)^k dv\\
&=&\frac{1}{\Gamma (\alpha)} \sum_{j=0}^k (-1)^{k-j} \left(\begin{array}{c} k\\ j \end{array} \right) \frac{\Gamma (\alpha +k -j)}{c^{\alpha +k -j}} \log^jx.
\end{eqnarray*}
Putting
$$B_\alpha (k,j) := \frac{\Gamma (\alpha +k-j)}{\Gamma (\alpha)} = \prod_{\nu =1}^{k-j} (\alpha +k -j -\nu),$$
we finally obtain
$$(J^\alpha_{0+,c}g_k)(x) = \sum_{j=0}^k (-1)^{k-j} \left(\begin{array}{c} k\\ j \end{array} \right)
\frac{B_\alpha(k,j)}{c^{\alpha +k -j}} \log^jx.$$

For the fractional derivative, putting $m = [\alpha] +1,$ we have
$$(D^\alpha_{0+,c} g_k)(x) = x^{-c}\sum_{j=0}^k (-1)^{k-j} \left(\begin{array}{c} k\\ j \end{array} \right)\frac{B_{m-\alpha}(k,j)}{c^{m-\alpha +k -j}}
\delta^m(x^c \log^j x).$$
In particular for $k=1,$ we have
\begin{eqnarray*}
(D^\alpha_{0+,c}g_1)(x) &=& x^{-c}\bigg[\frac{-B_{m-\alpha}(1,0)}{c^{m-\alpha+1}}\delta^m x^c + \frac{B_{m-\alpha}(1,1}{c^{m-\alpha}}
\delta^m (x^c \log x)\bigg]\\&=& x^{-c}\bigg[\frac{m-\alpha}{c^{m-\alpha+1}}\delta^m x^c + \frac{1}{c^{m-\alpha}}
\delta^m (x^c \log x)\bigg].
\end{eqnarray*}
Now using an easy induction, $\delta^m (x^c \log x) = c^m x^c \log x + mc^{m-1}x^c,$ thus we finally obtain the formula:
$$(D^\alpha_{0+,c}g_1)(x) = \alpha c^{\alpha -1}+ c^\alpha \log x.$$
This is another explanation of the fundamental theorem of fractional calculus in the Mellin frame. Indeed, it is easy to see that
$$J^\alpha_{0+,c}(D^\alpha_{0+,c}g_1)(x) = \log x.$$
We can obviously obtain formulae for higher values of $k.$

Note that the assumption $c >0$ is essential. Indeed as we remarked earlier, for $c=0,$ the function $\log x$ does not belong to the domain of the operator
$J^\alpha_{0+}.$

In \cite{KST}, Property 2.24, some related examples are treated concerning the Hadamard integrals $J^\alpha_{a+}f,$ with $a>0$ in place of $0.$

\item[3.] Let us consider the function $g(x) = e^{bt},$ $b \in \mathbb{R}.$ Then for any $c>0$ and $\alpha >0,$ we have $g \in Dom J^\alpha_{0+,c}$ and using the representation formula proved in \cite{BKT3}, Lemma 5(i), we have
$$(J^\alpha_{0+,c}g)(x) = \sum_{k=0}^\infty (c +k)^{-\alpha}\frac{b^k}{k!}x^k, \quad x \in \mathbb{R}^+$$
and the corresponding formula for the derivative, given in (Lemma 5(ii), \cite{BKT3}).
$$(D^\alpha_{0+,c}g)(x) = \sum_{k=0}^\infty (c +k)^{\alpha}\frac{b^k}{k!}x^k, \quad x \in \mathbb{R}^+.$$
As already remarked, assumption $c>0$ is essential. Indeed for $c=0$ Propositions \ref{analitica} and \ref{analitica 2} imply that we have similar representations for $J^\alpha_{0+}$ and $D^\alpha_{0+},$ only if the analytic function $f$ satisfies $f(0) = 0.$

Alternative representations are given by Propositions \ref{analitica 3}, \ref{analitica 4}, in terms of Stirling functions.
We have, for $\alpha >0,$
$$(D^\alpha_{0+,c}e^{bt})(x) = e^{bx}\sum_{k=0}^\infty S_c(\alpha, k) x^k b^k~~~~(x>0)$$
and
$$(J^\alpha_{0+,c}e^{bt})(x) = e^{bx}\sum_{k=0}^\infty S_c(-\alpha, k) x^k b^k~~~~(x>0).$$
\item[4.] Let us consider the "sinc" function which is analytic over the entire real line. The Taylor series is given by:
$$\mbox{sinc} (x) = \frac{\sin \pi x}{\pi x} = \sum_{k=0}^\infty (-1)^k\frac{\pi^{2k}}{(2k+1)!}x^{2k}.$$
Moreover it is easy to see that sinc $ \in X_{c, loc}$ for $c>0,$ while sinc $\not \in X_{0, loc}.$
Using Lemma 5 in \cite{BKT3} we have immediately
$$(J^\alpha_{0+,c}\mbox{sinc})(x) = \sum_{k=0}^\infty (-1)^k (c+2k)^{-\alpha} \frac{\pi^{2k}}{(2k+1)!}x^{2k}$$
and
$$(D^\alpha_{0+,c}\mbox{sinc})(x) = \sum_{k=0}^\infty (-1)^k (c+2k)^\alpha \frac{\pi^{2k}}{(2k+1)!}x^{2k}.$$
Another representation in term of Stirling functions of second type is a consequence of Proposition 13.

A formula for the (classical) derivatives of sinc can be found in \cite{BSS}.
For a given $s \in I\!\!N,$ differentiating by series, we have
\begin{eqnarray*}
(\mbox{sinc}~x)^{(s)} &=& \sum_{k=s}^\infty (-1)^k \frac{\pi^{2k}}{(2k+1)!}\frac{d^s}{dx^s}x^{2k} =
\sum_{k=s}^\infty  (-1)^k \frac{\pi^{2k}}{(2k+1)!}A_{s,k}x^{2k-s}\\
&=& \sum_{k=0}^\infty  (-1)^{k+s} \frac{\pi^{2k +2s}}{(2k+2s+1)!}A_{s,k+s}x^{2k+s},
\end{eqnarray*}
where
$$A_{s,k} = \prod_{\nu = 0}^{s-1} (2k -\nu).$$
Thus, using again Lemma 5 in \cite{BKT3}, we have
$$(J^\alpha_{0+,c}(\mbox{sinc}~t)^{(s)})(x) = \sum_{k=0}^\infty (-1)^{k+s}(c + 2k +s)^{-\alpha} \frac{\pi^{2k + 2s}}{(2k+2s +1)!}A_{s, k+s}x^{2k+s}$$
and
$$(D^\alpha_{0+,c}(\mbox{sinc}~t)^{(s)})(x) = \sum_{k=0}^\infty (-1)^{k+s}(c + 2k +s)^\alpha \frac{\pi^{2k + 2s}}{(2k+2s +1)!}A_{s, k+s}x^{2k+s}.$$
Note that for every odd $s,$ the above formula is valid also for $c=0,$ since in this instance $(\mbox{sinc}~x)^{(s)} \in X_{0,loc}.$
\end{enumerate}

\section{Applications to partial differential equations}

In this section we apply our theory to certain fractional differential equations. We notice here that the use of Mellin analysis in the theory of differential equations was considered in \cite{AB}, dealing with Cauchy problems for ordinary differential equations, involving Mellin derivatives of integral order. In \cite{EG}, Mellin analysis was applied to numerical solutions of Mellin integral equations. In the fractional case, differential equations were treated using various types of fractional derivatives, e.g. Riemann-Liouville, Caputo, Hadamard, etc (see \cite{KST}). The use of integral transforms is a very useful and used method for certain Cauchy or boundary value problems. However, the use of Mellin transforms in fractional differential equations involving Hadamard derivatives is so far not common.

Here we will examine certain boundary value problems related to an evolution equation and to a diffusion problem, using the Mellin transform approach and using Hadamard derivatives. In the first example, the fractional evolution equation originates from a Volterra integral equation with a special kernel. The second example is a fractional diffusion equation.

\subsection{An integro-differential equation}

Let $\alpha \in ]0,1[$ be fixed, and let
$$K_\alpha(x,u):=\frac{1}{\Gamma (1-\alpha)}\bigg(\log \frac{x}{u}\bigg)^{-\alpha}\chi_{]0,x[}(u),\quad x >0.$$
Let us consider the following problem: find a function $w:\mathbb{R}^+ \times \mathbb{R}^+ \rightarrow \mathbb{C}$ such that $w(x,0) = f(x),$ for a fixed boundary data $f:\mathbb{R}^+ \rightarrow \mathbb{C},$ and
\begin{eqnarray}
Ax \frac{\partial}{\partial x}\int_0^\infty K_\alpha(x,u) w(x,y)\frac{du}{u} + B \frac{\partial}{\partial y}w(x,y) = 0,
\end{eqnarray}
$A, B$ being two positive constants.

Now, equation (14) can be rewritten as a fractional partial differential evolution equation in the Hadamard sense, as
$$A (D^\alpha_{0+}w(\cdot, y))(x) + B \frac{\partial}{\partial y}w(x,y) = 0,~~~(x,y \in \mathbb{R}^+).$$

Without loss of generality we can assume $A = B = 1,$ thus
\begin{eqnarray}
(D^{\alpha}_{0+} w(\cdot, y))(x)= -\frac{\partial}{\partial y} w(x,y),~~(x,y \in \mathbb{R}^+)
\end{eqnarray}
with initial data $w(x,0) = f(x),~x>0.$
We call for a function $w: \mathbb{R}^+ \times \mathbb{R}^+ \rightarrow \mathbb{C}$ satisfying the following properties
\begin{description}
\item[1)]
$w(\cdot,y) \in {\mathcal X}^{\alpha}_{[a, 0]}$ for every $y>0$ and for a fixed $a<0$
\item[2)]
there is a function $K\in X_\nu,$ $\nu \in [a,0[,$ such that for every $x,y>0$
$$ \bigg| \frac{\partial }{\partial y} w(x,y) \bigg| \leq K(x) $$
\item[3)] for a fixed $f \in X_\nu,$ we have
$ \lim_{y\rightarrow 0^+} \| w(\cdot, y) - f(\cdot) \|_{X_\nu} =0.$
\end{description}
Assuming that such a function exists we apply the Mellin transform with respect to the variable $x$ on the line $\nu + it$ to both sides of (15), obtaining
$$ [D^{\alpha}_{0+} w(\cdot, y)]^{\wedge}_M(\nu +it)= -\bigg[\frac{\partial}{\partial y} w(\cdot, y)\bigg]^{\wedge}_M(\nu +it).$$
Using Theorems 1 and 5 we have
$$ [D^{\alpha}_{0+} w(\cdot, y)]^{\wedge}_M(\nu +it)=(-\nu -it)^{\alpha}[w(\cdot, y)]^{\wedge}_M(\nu +it).$$
Moreover by property 2),
$$  \bigg[\frac{\partial}{\partial y} w(\cdot, y)\bigg]^{\wedge}_M(\nu +it)=\int_0^{+\infty} x^{\nu +it-1} \frac{\partial}{\partial y} w(x,y) dx =  \frac{\partial}{\partial y} [w(\cdot,y)]^{\wedge}_M (\nu +it),$$
thus equation (15) is transformed into a first order ordinary differential equation
$$ (-\nu -it)^{\alpha}[w(\cdot,y)]^{\wedge}_M(\nu +it) =  -\frac{\partial}{\partial y} [w(\cdot,y)]^{\wedge}_M(\nu +it)$$
which has the solution
$$ [w(\cdot,y)]^{\wedge}_M(\nu +it) = A(\nu +it) e^{-(-\nu -it)^{\alpha}y}$$
where $A(\nu +it)$ is independent of $y.$
The determination of $A(\nu +it)$ follows from condition 3); indeed we have that
$ [w(\cdot,y)]^{\wedge}_M(\nu +it) \rightarrow [f]^{\wedge}_M (\nu +it)$ uniformly for $y\rightarrow 0^+$ and for $t\in \mathbb{R},$ and so $A(\nu +it) = [f]^{\wedge}_M (\nu +it),$ obtaining
$$  [w(\cdot,y)]^{\wedge}_M(\nu +it)= [f]^{\wedge}_M (\nu +it) e^{-(-\nu-it)^{\alpha}y}.$$
Now putting $s = -\nu -it,$ we have $Re s = -\nu >0$ and so, since $y>0,$ the inverse Mellin transform of $e^{-y s^{\alpha}}$ exists and it is given by (see Theorem 6 in \cite{BJ})
\begin{eqnarray}
G(x,y):= \frac{1}{2\pi}\int_{-\infty}^{+\infty} e^{-(-\nu -it)^\alpha y}x^{-\nu -it}dt.
\end{eqnarray}
Thus if the solution of (15) exists, by the Mellin-Parseval formula (see \cite{BJ}), it has the form
\begin{eqnarray*}
w(x,y) = \int_0^{+\infty} f(v) G(\frac{x}{v},y)\frac{dv}{v},~~x,y>0.
\end{eqnarray*}
In order to verify that the function $w(x,y)$ is actually a solution of the problem we make a direct substitution.
We have, by differentiating under the integral
\begin{eqnarray*}
&&-\frac{\partial w}{\partial y}(x,y) = \int_0^{+\infty} f(v) \bigg[\frac{1}{2\pi}\int_{-\infty}^{+\infty}\bigg(\frac{x}{v}\bigg)^{-\nu -it}(-\nu -it)^\alpha e^{-(-\nu -it)^\alpha y}dt\bigg]\frac{dv}{v}\\ &=&
\frac{1}{2\pi}\int_{-\infty}^{+\infty}(-\nu -it)^\alpha x^{-\nu -it} e^{-(-\nu -it)^\alpha y}[f]^\wedge_M(\nu +it)dt.
\end{eqnarray*}
Now, let us consider
$$(D^{\alpha}_{0+} w(\cdot, y))(x) = \delta (J^{1-\alpha}_{0+}w(\cdot,y))(x).$$
We have
\begin{eqnarray*}
&&(D^{\alpha}_{0+} w(\cdot, y))(x) = (x \frac{\partial}{\partial x})\bigg[\frac{1}{\Gamma(1-\alpha)} \int_0^x
\bigg(\log \frac{x}{u}\bigg)^{-\alpha}\bigg(\int_0^{\infty} f(v) G(\frac{u}{v},y)\frac{dv}{v}\bigg)\frac{du}{u}\bigg]\\
&=&(x \frac{\partial}{\partial x})\bigg[\frac{1}{\Gamma(1-\alpha)} \int_1^{+\infty}
(\log z)^{-\alpha}\bigg(\int_0^{\infty} f(v) G(\frac{x}{zv},y)\frac{dv}{v}\bigg)\frac{dz}{z}\bigg]\\&=&
\frac{x}{\Gamma(1-\alpha)} \int_1^{+\infty}
(\log z)^{-\alpha}\bigg(\int_0^{\infty} f(v) \frac{\partial}{\partial x}G(\frac{x}{zv},y)\frac{dv}{v}\bigg)\frac{dz}{z}.
\end{eqnarray*}
Since
$$\frac{\partial}{\partial x}G(x,y) = \frac{1}{2\pi}\int_{-\infty}^{+ \infty}e^{-(-\nu -it)^\alpha y}(-\nu -it)x^{-\nu -it -1}dt,$$
 putting $s = -\nu -it,$ we obtain
\begin{eqnarray*}
&&(D^{\alpha}_{0+} w(\cdot, y))(x)\\&=& \frac{x}{\Gamma(1-\alpha)} \int_1^{+\infty}
(\log z)^{-\alpha}\bigg(\int_0^{\infty} f(v) \frac{1}{zv} \bigg(\frac{1}{2\pi}\int_{-\infty}^{+ \infty}e^{-s^\alpha y}s\bigg(\frac{x}{zv}\bigg)^{s -1}dt\bigg)\frac{dv}{v}\bigg)\frac{dz}{z}\\
&=&\frac{1}{2\pi}\frac{1}{\Gamma(1-\alpha)} \int_1^{+\infty}
(\log z)^{-\alpha}\bigg(\int_{-\infty}^{+\infty} e^{-s^\alpha y}s \bigg(\frac{x}{z}\bigg)^{s}[f]^\wedge_M(-s)dt\bigg) \frac{dz}{z} \\
&=&
\frac{1}{2\pi}\frac{1}{\Gamma(1-\alpha)} \int_{-\infty}^{+\infty} [f]^\wedge_M(-s) e^{-s^\alpha y}s\bigg(\int_0^{x}
\bigg(\log \frac{x}{u}\bigg)^{-\alpha} u^{s} \frac{du}{u}\bigg)dt.
\end{eqnarray*}
Since Example 1 of Section 6, holds for $c=0$ and complex $b$ with Re $ b >0,$ we have
$$\frac{1}{\Gamma (1-\alpha)} \int_0^x \bigg(\log \frac{x}{u}\bigg)^{-\alpha} u^s \frac{du}{u} = (J^{1-\alpha}_{0+}u^s)(x) =
s^{\alpha -1}x^s,$$
and so we have
$$(D^{\alpha}_{0+} w(\cdot, y))(x) = \frac{1}{2\pi}\int_{-\infty}^{+\infty}s^\alpha x^{s} e^{-s^\alpha y}[f]^\wedge_M(\nu +it)dt,$$
i.e. the assertion. So we have proved the following
\begin{Theorem}\label{equadiff}
Under the assumptions imposed, equation (15) with the initial data $f,$ has the unique solution given by
$$w(x,y) = \int_0^{+\infty} f(v) G(\frac{x}{v},y)\frac{dv}{v},~~x,y>0,$$
where the function $G(x,y)$ is defined in (16).
\end{Theorem}
\vskip0,3cm
Note that for $\alpha = 1/2,$ we have a closed form for the function $G(x,y).$ Indeed, using
formula 3.7 page 174 in \cite{OB}, we obtain
$$G(x,y) = \frac{y}{2\sqrt{\pi}} (-\log x)^{-3/2}\mbox{exp}\bigg(\frac{y^2}{4 \log x}\bigg)\chi_{]0,1[}(x)$$
and the solution is then given by
$$w(x,y) = \frac{y}{2\sqrt{\pi}} \int_0^{1} f(\frac{x}{v}) (-\log v)^{-3/2}\mbox{exp}\bigg(\frac{y^2}{4 \log v}\bigg)\frac{dv}{v}.$$
Equation (15) was also discussed in \cite{KST}, but using fractional Caputo derivatives. Our treatment, however, contains real proofs.


\subsection{A diffusion equation}

For $\alpha >0$ let us consider the fractional diffusion equation
\begin{eqnarray}
(D^\alpha_{0+} w(\cdot, y))(x) = \frac{\partial^2}{\partial y^2}w(x,y), ~~~(x,y \in \mathbb{R}^+)
\end{eqnarray}
with the initial condition
$$\lim_{y \rightarrow 0^+}\|w(\cdot,y) - f(\cdot)\|_{X_0} = 0,$$
for a fixed $f \in X_0.$

We call for a function $w: \mathbb{R}^+\times \mathbb{R}^+ \rightarrow \mathbb{C}$ satisfying the following assumptions:
\begin{description}
\item[1)]
$w(\cdot,y) \in {\mathcal X}^{\alpha}_{0}$ for every $y>0,$ and there exist $N>0$ such that $\|w(\cdot,y)\|_{X_0} \leq N,$ for every
$y \in \mathbb{R}^+.$
\item[2)]
there are functions $K_1, K_2 \in X_0,$ such that for every $x,y>0$
$$ \bigg| \frac{\partial }{\partial y} w(x,y) \bigg| \leq K_1(x),~ \bigg| \frac{\partial^2 }{\partial y^2} w(x,y) \bigg| \leq K_2(x)$$
\item[3)] for a fixed $f \in X_0,$ we have
$ \lim_{y\rightarrow 0^+} \| w(\cdot, y) - f(\cdot) \|_{X_0} =0.$
\end{description}
Using the same approach as in the previous example, taking the Mellin transforms of both sides of the equation (17), we obtain
$$ [D^{\alpha}_{0+} w(\cdot, y)]^{\wedge}_M(it)= -\bigg[\frac{\partial^2}{\partial y^2} w(\cdot, y)\bigg]^{\wedge}_M(it).$$
Using Theorems 1 and 5 we have
$$ [D^{\alpha}_{0+} w(\cdot, y)]^{\wedge}_M(\nu +it)=(-it)^{\alpha}[w(\cdot, y)]^{\wedge}_M(it).$$
Moreover by property 2),
$$\bigg[\frac{\partial^2}{\partial y^2} w(\cdot, y)\bigg]^{\wedge}_M(it)= \frac{\partial^2}{\partial y^2} [w(\cdot,y)]^{\wedge}_M (it),$$
thus equation (17) is transformed into the second order linear ordinary differential equation
\begin{eqnarray}
(-it)^\alpha z_t(y) = z''_t(y),\quad y >0,
\end{eqnarray}
with respect to the function
$$z_t(y):= [w(\cdot,y)]^\wedge_M(it),~~~t \in \mathbb{R}.$$
If $t=0$ the solution is the linear function $z_0(y) = A(0) + B(0)y,$ while for $t \neq 0,$ the characteristic equation
associated with (18)
$$\lambda^2 = \exp (\alpha \log (-it)),$$
has two complex solutions
$$\lambda_1:= |t|^{\alpha/2}\bigg(\cos \frac{\alpha \pi}{4} +i \sin \frac{\alpha \pi}{4}(-\sgn t)\bigg),$$
$$\lambda_2:= -|t|^{\alpha/2}\bigg(\cos \frac{\alpha \pi}{4} +i \sin \frac{\alpha \pi}{4}(-\sgn t)\bigg).$$
Thus, for $t \neq 0,$ we obtain the general solution
$$z_t(y) = A(t) e^{-|t|^{\alpha/2}(\cos (\alpha \pi/4) + i(-\sgn t)\sin (\alpha \pi/4))y} +
B(t) e^{|t|^{\alpha/2}(\cos (\alpha \pi/4) + i(-\sgn t)\sin (\alpha \pi/4))y}.$$
Now, let $\alpha$ be such that $\cos (\alpha\pi/4) >0.$
By the boundary condition 3), we have also that $z_t(y)$ is uniformly convergent to $[f]^\wedge_M$ as $y\rightarrow 0^+.$
Moreover, by assumption 1), there exists a constant $N>0$ such that
$|z_t(y)| \leq N,$ for every $t \in \mathbb{R}.$ This means that we must have $B(t) = 0$ for
every $t \in \mathbb{R},$ thus
$$z_t(y) = [w(\cdot, y)]^\wedge_M(it) = [f]^\wedge_M(it) e^{-|t|^{\alpha/2}(\cos (\alpha \pi/4) + i(-\sgn t)\sin (\alpha \pi/4))y}.$$
Now, the function
 $$e^{-|t|^{\alpha/2}(\cos (\alpha \pi/4) + i(-\sgn t)\sin (\alpha \pi/4))y}$$
is summable as a function of $t \in \mathbb{R},$ and its inverse Mellin transform is given by
$$G(x,y):= \frac{1}{2\pi}\int_{-\infty}^\infty e^{-|t|^{\alpha/2}(\cos (\alpha \pi/4) + i(-\sgn t)\sin (\alpha \pi/4))y} x^{-it}dt.$$
Then if a solution exists it has the form
$$w(x,y) = \int_0^\infty f(u) G(\frac{x}{u}, y)\frac{du}{u},~~~~x,y >0.$$
Analogously, if $\alpha$ is such that $\cos (\alpha \pi/4) <0,$ then we have $A(t)=0$ for every $t \in \mathbb{R},$ and the corresponding function
$G(x,y)$ takes the form
$$G(x,y) = \frac{1}{2\pi}\int_{-\infty}^\infty e^{-|t|^{\alpha/2}(|\cos (\alpha \pi/4)| + i(-\sgn t)\sin (\alpha \pi/4))y} x^{-it}dt.$$
That the above function is really a solution can be proved, as before, by a direct subsitution into the differential equation.

The function $G(x,y)$ can be written in a more simple form. Indeed, using Euler's formula, putting $a:= |\cos (\alpha \pi/4)|,$ $b:= \sin (\alpha \pi/4),$
we can write:
\begin{eqnarray*}
&&G(x,y) = \frac{1}{2\pi}\int_0^\infty e^{-|t|^{\alpha/2}(a-ib)y}(\cos (t \log x - i\sin (t \log x))dt \\
&+&
\frac{1}{2\pi}\int_0^\infty e^{-t^{\alpha/2}(a+ib)y}(\cos (t \log x + i\sin (t \log x))dt\\
&=& \frac{1}{\pi}\int_0^\infty e^{-t^{\alpha/2}ay}[\cos (t \log x)\cos (t^{\alpha/2}by + \sin (t\log x) \sin (t^{\alpha/2}by)]dt\\
&=&\frac{1}{\pi}\int_0^\infty e^{-t^{\alpha/2}ay}\cos (t \log x - t^{\alpha/2}by)dt.
\end{eqnarray*}

For $\alpha = 1$ using Proposition 13, we obtain the (not fractional) equation
$$x\frac{\partial w}{\partial x}(x,y) = \frac{\partial^2 w}{\partial y^2}(x,y),~~~~(x,y \in \mathbb{R}^+$$
and using our approach the corresponding problem has a unique solution of the form
$$w(x,y) = \int_0^\infty f(u) G(\frac{x}{u}, y)\frac{du}{u},~~~~x,y >0,$$
where
$$G(x,y) = \frac{1}{\pi}\int_0^\infty \exp(-\frac{\sqrt{2t}y}{2})\cos (\frac{\sqrt{2t}y}{2}- t \log x) dt.$$
This integral has a closed form. Indeed by an elementary subsitution, we can write
$$I:=\int_0^\infty \exp(-\frac{\sqrt{2t}y}{2})\cos (\frac{\sqrt{2t}y}{2}- t \log x) dt =
2\int_0^\infty \exp(-\frac{\sqrt{2}yu}{2})\cos (\frac{\sqrt{2}yu}{2}- u^2 \log x) dt.$$
Now, the above integral, depending on the sign of $\log x,$ can be reduced to the integrals ($p>0)$ (see \cite{GR}, page 499)
$$\int_0^\infty v e^{-pv}\cos (2v^2 -pv)dv = \frac{p\sqrt{\pi}}{8}\exp(-p^2/4),$$
if $\log x>0$ and
$$\int_0^\infty v e^{-pv}\cos (2v^2 + pv)dv = 0,$$
if $\log x \leq 0.$
Indeed, if we put $u= \sqrt{2/\log x}v$ in the first case, and $u=\sqrt{2/|\log x|}v$ in the second case, we get easily
$$I = \frac{\sqrt{\pi}}{\log x\sqrt{2\log x}}\exp \bigg(-\frac{y}{2\sqrt{2}\log x}\bigg),~~~x>1,$$
while $I=0$ for $0<x<1.$ Therefore,
\begin{eqnarray*}
G(x,y) = \left\{ \begin{array}{ll} \sqrt{\displaystyle\frac{\pi}{2}}\displaystyle\frac{1}{(\log x)^{3/2}}\exp\bigg(-\displaystyle\frac{y}{2\sqrt{2}\log x}\bigg), & ~x>1,~y>0 \\ \\
0, & ~0<x\leq 1,~y>0.
\end{array} \right.
\end{eqnarray*}
For $\alpha = 1/2,$ the equation becomes
$$(D^{1/2}_{0+} w(\cdot, y))(x) = \frac{\partial^2}{\partial y^2}w(x,y),~~~~(x,y \in \mathbb{R}^+)$$
and putting
$a:= \cos (\pi/8) = \sqrt{2 + \sqrt{2}}/2, ~b:= \sin (\pi/8) = \sqrt{2 - \sqrt{2}}/2,$ we obtain the following representation of the function $G(x,y):$
$$G(x,y) = \frac{1}{\pi}\int_0^\infty \exp(-\sqrt[4]{t}ay) \cos (t \log x - \sqrt[4]{t}by)dt.$$

For $\alpha = 4,$
using Proposition 13, our equation equation  has the form
\begin{eqnarray}
\sum_{k=0}^4 S_0(4,k) x^k\bigg(\frac{\partial}{\partial x}\bigg)^{(k)}w(x,y) = \frac{\partial^2 w}{\partial y^2}(x,y),~~(x,y \in \mathbb{R}^+)
\end{eqnarray}
i.e.
$$x^4\frac{\partial^4 w}{\partial x^4}(x,y)+ 6x^3\frac{\partial^3w}{\partial x^3}(x,y) +
7x^2\frac{\partial^2 w}{\partial x}(x,y) +
x\frac{\partial w}{\partial x}(x,y) =\frac{\partial^2 w}{\partial y^2}(x,y),~~(x,y \in \mathbb{R}^+)$$
In this instance we have $\cos (\alpha\pi/4) = -1,$ and so,
the unique solution of our problem for equation (19) has the form
$$w(x,y) = \int_0^\infty f(u) G(\frac{x}{u}, y)\frac{du}{u},~~~~x,y >0,$$
where
$$G(x,y) = \frac{1}{\pi}\int_0^\infty e^{-t^2y}\cos (t \log x) dt.$$
This integral can be reduced by an elementary substitution, to the classical integral
$$g(v) = \int_0^\infty e^{-t^2} \cos (tv) dt = \frac{\sqrt{\pi}}{2}\exp(-v^2/4),$$
thus obtaining
$$G(x,y) = \frac{1}{2}\sqrt{\frac{\pi}{y}}\exp(-\log^2 x/4y).$$
Another example in fractional case, is $\alpha = 5/2.$  In this case we have
$a = |\cos ((5/8)\pi)|= \sqrt{2 - \sqrt{2}}/2$ and $b= \sin ((5/8)\pi) = \sqrt{2 + \sqrt{2}}/2.$ The corresponding
function $G(x,y)$ is given by
$$G(x,y) = \frac{1}{\pi}\int_0^\infty \exp({-t^{5/8}\sqrt{2-\sqrt{2}}y})\cos (t \log x - t^{5/4}\sqrt{2+\sqrt{2}}y/2)dt.$$
The above approach works for every value of $\alpha$ except those for which $\cos(\alpha \pi/4) = 0.$ For $\alpha = 2,$
the resulting wave equation in the Mellin setting reads
$$x^2 \frac{\partial^2}{\partial x^2}w(x,y) + x \frac{\partial}{\partial x}w(x,y) = \frac{\partial^2}{\partial y^2}w(x,y),~~(x,y \in \mathbb{R}^+)$$
But this equation is treated in detail in \cite{BJ} with different boundary conditions. Experts in the evaluations of integrals could surely obtain more elegant representations of the $G(x,y)-$ functions.


\section{ A short biography of R.G. Mamedov and some historical notes}

Rashid Gamid-oglu Mamedov (changed into Mammadov since 1991), born in a peasant family on December 27, 1931, in the village Dashsalakhly, Azerbaijan SSR, lost
his father at the age 6 and grew up with his mother and three sisters.

After spending the school years 1938-48 in the middle school of his home village, he was admitted to the Azerbaijan Pedagogical Institute (API) in Baku.
In 1952, he graduated from its Mathematics Department with a so-called red diploma-honours (i.e. diploma cum laude). Immediately he was accepted for post-graduate study
at the Chair of Mathematical Analysis of API, and defended his PhD thesis ("Kandidatskaya") entitled "Some questions of approximation by entire functions and polynomials" in 1955. This dissertation was one basis to the monograph "Extremal Properties of Entire Functions" published in 1962 by his scientific supervisor I.I. Ibragimov. During the years 1953-1960, R.G. Mamedov was affiliated with the Chair of Mathematical Analysis at API in various positions, first as assistant (1953-1956) and senior lecturer (1956-57), later as docent (assistant professor, 1957-1960)

In 1960-1963, R.G. Mamedov held a position as senior researcher at the Institute of Mathematics and Mechanics of the Azerbaijan Academy of Science. Free of teaching duties, he published in a very short period of time his fundamental contributions to the theory of approximation by linear operators which made him known both in the former Soviet Union and abroad. These deep results comprised his "Doktorskaya" (Habilitation degree) "Some questions of approximation of functions by linear operators" submitted to Leningrad State Pedagogical A.I.Herzen-Institute in 1964. At the age of 33 years, R.G. Mamedov was awarded the Dr. of Phys. and Math. degree and was appointed as full professor to the Chair of Higher Mathematics at Azerbaijan Polytechnic Institute in Baku. Here he started his remarkable career as university teacher and educator, supervising as many as 23 PhD theses over the years, two of his students obtained the Dr. of Phys. and Math. degree themselves.
In 1966, he gave a contributed talk at the ICM Congress in Moscow. \\

In 1967, he published his first monograph "Approximation of Functions by Linear Operators",
recognised by the international mathematical community, although it was written in Azerbaijani.
His son Aykhan reported that his father possessed a copy of \cite{BN} and recalls him speaking about the authors.
In 1969, R.G. Mamedov, was appointed head of the Chair of Higher Mathematics at Azerbaijan State Oil Academy in Baku, a position which he held for 26 years. His cycle of investigations on properties of integral transforms of Mellin-type led to the publication of several research monographs, in particular "On Approximation of Conjugate Functions by Conjugate M-Singular Integrals" (1977), "On Approximation of Functions by Singular Integrals of Mellin Type" (1979), and
"Mellin Transform and Theory of Approximation" (1991). With equal enthusiasm, he created textbooks for use at the Azerbaijan institutions of higher education that are still of widespread use. His three-volume "Course of Higher Mathematics" (1978, 1981, 1984) has several editions. R.G. Mamedov is also the author
of 20 booklets and articles popularising mathematics among the general public and raising the standards of mathematics education in his home country.

R.G. Mamedov was not only an outstanding scientist and educator but also impressed everybody who met him by his outgoing character, friendly personality,
and for being very accessible and supportive in personal and scientific matters. He married in 1960, two of his three sons being mathematicians themselves.
R.G. Mamedov died on May 2, 2000, at the age of 68 after an infarct. He is survived by his 
spouse Flora Mamadova and three sons, there now being  seven grandchildren, five boys and two girls, four being born after his death.

 \begin{figure}[!h]
                \centering
               \fbox{ \includegraphics[width=0.6\textwidth]{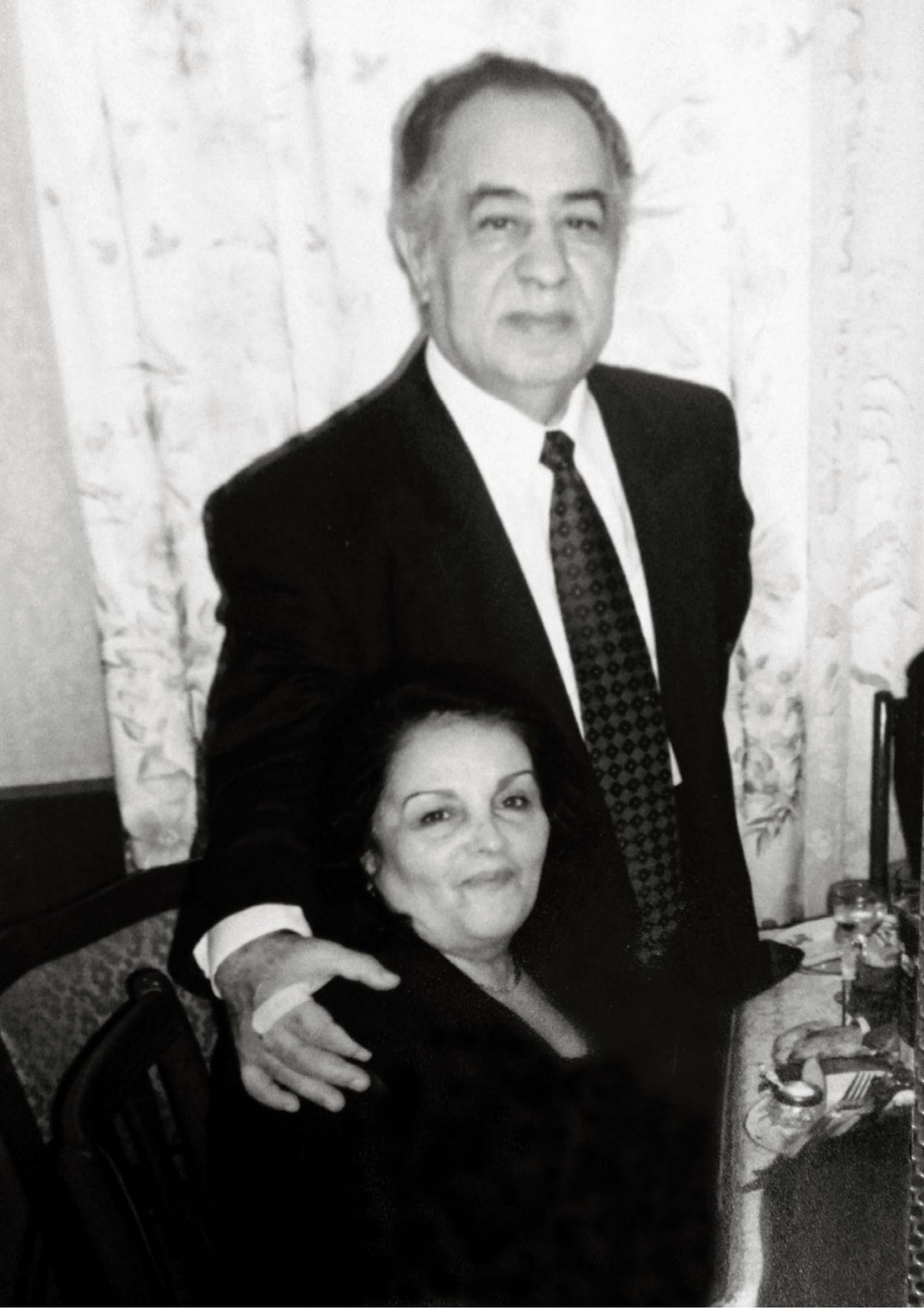}}
                \caption{\small  A photo of Prof.  Rashid Mamedov together with his spouse Flora Mamedova, who now takes her husband's role in keeping alive Azerbaijani customs among her grandchildren. It was taken in the year of his death, 2000.}
                \label{nome etichetta se vuoi richiamare la foto}
        \end{figure}

Work in the broad area of approximation theory at the University of Perugia, was initiated by its former visionary, departmental director,
 C. Vinti (1926-1997) a master in Calculus of Variation (see \cite{VIN}). It was decisively influenced by the work of J. Musielak, a chief representative of the Orlicz analysis school at Poznan, its first joint work bring in the direction of (nonlinear) integral operators in the setting of modular spaces (\cite{BMV}), as well as by the work at Aachen, together with P. L. Butzer and R. L. Stens. During recent research at Perugia  in matters asymptotic expansions of certain Mellin-type convolution operators and convergence properties in the spaces of functions of bounded variation (see \cite{AV1,AV2}, \cite{BM, BM1, BM2, BM3, BM4, BM5, BM6}), a MathSciNet search led to the treatise of R.G. Mamedov under discussion. Since it was nowhere to be found, it was finally A. Gadjiev, Academy of Sciences of Azerbaijan, who within a few weeks kindly sent a copy, as a present.
It has served us well not only in our local work at Perugia but also in the present joint investigation.

As to the work at Aachen, although we knew of the existence of the great school of approximation theory at Leningrad since 1949 (through G.G.Lorentz),
it was the Second All-Union Conference on Constructive Theory of Functions, held at Baku on Oct.8-13, 1962 , that drew our attention
to approximation theory at Baku. That was a couple of years after its  proceedings (with 638 pp.) appeared in 1965.
(The Aachen group organised the
 first conference  on approximation in the West (August 4-10,1963 ;ISNM, Vol. 5, Birkhaeuser,
Basel, 1964)).

It was Aachen's  former student E.L. Stark
(1940-1984), who in view of his fluent knowledge of Russian kept well aware of approximation theoretical studies at Leningrad,
Moscow and Kiev, was surprised when he discovered the Baku proceedings. In fact, Russian approximation theory was a model for us in Aachen, especially in its earlier years; and Stark's great input benefited us all.
We exchanged letters with R.G.Mamedov and in 1974 invited him to participate in our Oberwolfach conference on Linear Operators and Approximation II,
held March 30 - April 6. But he was unable to attend at the last moment (likewise in the case of S.M. Nikolskii, S.A. Teljakovski and B.S. Mitijagin), as is recorded in its Proceedings (ISNM, Vol. 25, Birkhaeuser Basel, 1974).
In our volume with R.J Nessel," Fourier Analysis
and Approximation (Birkhaeuser/Academic Press, 1971), we cited eight papers
of R.G. Mamedov , plus his book "Approximation of Functions by Linear Operators" (Azerbaijani,
 Baku, 1966).  They played a specific role in our book. The work on Mellin analysis at Aachen , together with S. Jansche (see \cite{BJ}, \cite{BJ1}, \cite{BJ2}, \cite{BJ3}) was independent of that at Baku.

\thispagestyle{empty}
\begin{sidewaysfigure}
\centering
\fbox{\includegraphics[angle=90,width=0.85\textwidth]{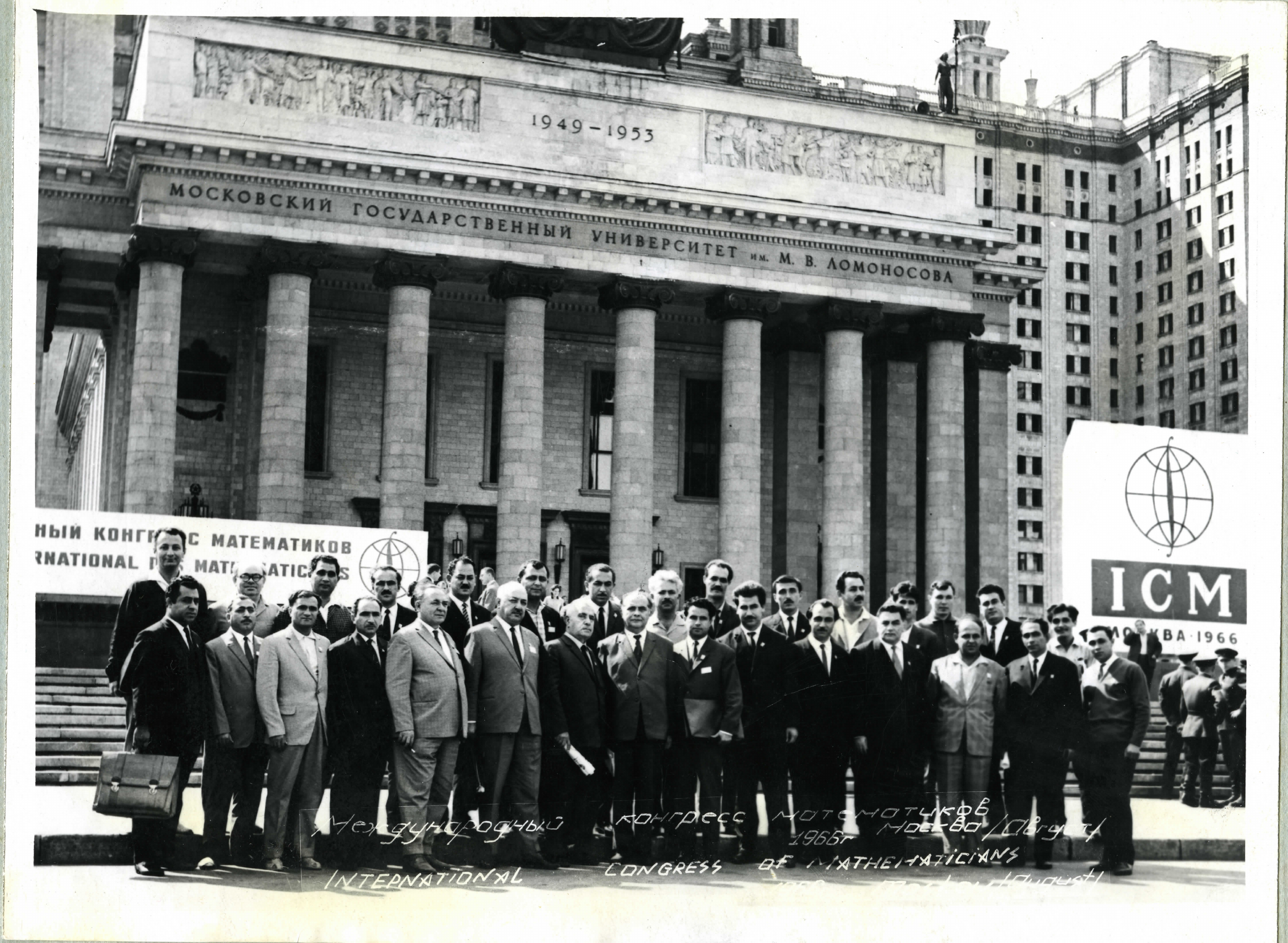}}
\caption{\small Private photo of the  30 Azerbaijani participants at  the  ICM, held in Moscow 1966, and kindly
	      forwarded to the authors by  Prof. Boris Golubov. Prof.  Mamedov stands with his large briefcase in the first row,
		 on the extreme left, the President of the Azerbaijani Academy of Sciences (in 1966), Acad,Z. Khalilov, stands in the      
         center of the first row, eighth from the left, together with  Prof.I.I. Ibragimov  (fifth  from the left) and the Dean of the Mechanical-    
         Mathematical Department of  Azerbaijani State University, Prof. A.I. Guseinov, sixth from the left. Prof.Golubov was invited to
         be present in this photo  since he spent his first three years  (1956-195 ) as a student  at their university,
          and also participated. in the Congress. We find him in the second row, third from the right.}\label{landscapefig}
\end{sidewaysfigure}
\newpage

\section{Concluding remarks}

The theory of Mellin analysis is a fascinating field of research, one still in the state of development, one which will surely have further important applications in various
fields of applied mathematics. As noted in the Introduction, a pioneering  contribution in this direction was the treatise of R.G. Mamedov \cite{MA}.
The translation into English of the main part of Mamedov's preface reads: {\em In classical approximation theory approximation of functions by polynomials and entire functions are considered, and relations between the order of best approximation of the functions and their structural and differential properties are studied. In connection with the saturation problem and P.P. Korovkin theorems on the convergence of linear positive operators, numerous investigations are dedicated to the approximations of functions by linear operators, in particular by linear positive operators, and by various singular integral operators.
To this aim some function classes are introduced and studied. Moreover, the saturation classes of different linear operators by means of Fourier transform or other integral transforms are investigated. Many results in this field and the base of the theory of integral Fourier transform were published in the fundamental monograph of P.L. Butzer and R.J. Nessel "Fourier analysis and approximation". At present some other integral transforms are also used in studying different function classes and the associated saturation order of approximation by linear operators.

The Mellin transform has important applications in the solution of boundary value problems in wedge shaped regions. It is also one of the most important methods for the study of classes of functions defined on the positive real line.
The theory of Mellin transform requires the introduction of new concepts of derivative and integral, called M-derivative and M-integral.
In this field in recent years many results have been produced.
In this monograph we attempt significantly to complement those 
results and introduce them from the unified point of view. I have used 
material written earlier in the book with G.N. Orudzhev, namely "On the approximation of functions by singular integrals of Mellin type, Baku, 1979.}

After that, Mellin analysis was introduced in a systematic way in \cite{BJ}, \cite{BJ1}, \cite{BJ2}, then developed in
\cite{BKT}, \cite{BKT1}, \cite{BKT2}, \cite{BKT3}, \cite{BKT4} and later on in \cite{BM}, \cite{BM1}, \cite{BM2}, \cite{BM3}, \cite{BM4},
\cite{BM5}, \cite{BM6}, \cite{MAN}.
Many other results and applications are surely to be discovered and the present paper is a further contribution in this direction.

Our theory of Hadamard-type fractional integrals and derivatives is concerned with real values of the parameter $\alpha.$ The extension to complex
values of $\alpha$ can be carried out essentially in the same way, assuming Re $\alpha \geq 0$ in place of $\alpha \geq 0$ (see also e.g. \cite{BKT3}). For general complex values $\alpha \in \mathbb{C},$ the theory may be more delicate. As an example, in Theorem \ref{Theorem 1}, the assumption Re $\alpha >0$ is basic for the application of the Abel-Stolz theorem. Indeed, for complex values of $\alpha$ such that Re $\alpha <0$ the convergence of the binomial series on the boundary of its convergence disk may fail.
For Re $\alpha \leq -1,$ this convergence fails at every point of the boundary, while for $-1< Re~ \alpha <0,$ it fails at just one point.
\vskip0,5cm
\noindent
{\bf Aknowledgments} The authors would like to thank Boris Ivanovich Golubov (Moscow) for his great help in regard to the short biography of R.G. Mamedov. He contacted his colleague in Baku, who in turn received a four pages biography of Mamedov together with a list of his publications kindly sent by his son Aykhan Mammadov, and made a first translation of this biography. The present biography is the extended and polished version kindly carried out by Peter Oswald (Bremen). 
The authors are grateful to Aykhan Mammadov for his extraordinary help in regard to various aspects
of his fathers life, Azerbaijan, and the paper itself.
The translation of the preface of \cite{MA} is due to Andi Kivinukk (Tallinn). Also, the authors wish to thank Annarita Sambucini for her technical support in including photos in the text.

The first and third authors have been partially supported by the Gruppo
Nazionale Analisi Matematica, Probabilità e Applicazioni (GNAMPA) of the Istituto Nazionale di Alta Matematica
(INdAM), through the INdAM - GNAMPA Project 2014,  and by the Department of Mathematics and Computer Sciences of University of Perugia.


\end{document}